\documentclass[12pt]{article}
\usepackage{bbm}
\usepackage{amsfonts}
\usepackage{pifont}
\usepackage{hyperref}
\usepackage{mathrsfs}
\usepackage{indentfirst}
\usepackage{amsmath}
\usepackage{amssymb}
\usepackage[amsmath, thmmarks]{ntheorem}
\usepackage{cite}
\usepackage{graphicx}
\usepackage{tikz}
\newtheorem{definition}{\bf Definition}[section]

\newtheorem{theorem}{\bf Theorem}[section]
\newtheorem{remark}{\bf Remark}[section]

\newtheorem{proposition}{\bf Proposition}[section]
\newtheorem{example}{\bf Example}[section]

\begin{document}
\setcounter{page}{1}

\title{{\textbf{New construction methods for uninorms via functions with $q$ and uninorms on bounded lattices}}\thanks {Supported by National Natural Science
Foundation of China (No.11871097, 12271036)}}
\author{Zhen-Yu Xiu$^1$\footnote{\emph{E-mail address}: xyz198202@163.com },  Zheng-Yuan Si$^1$\footnote{\emph{E-mail address}: 1161559146@qq.com }, Xu Zheng$^3$\footnote{Corresponding author. \emph{E-mail address}: 3026217474@qq.com }\\
\emph{\small $^{1,2,3}$  College of Applied Mathematics, Chengdu University of Information Technology, }\\
\emph{\small Chengdu 610000,  China }}
\newcommand{\pp}[2]{\frac{\partial #1}{\partial #2}}
\date{}
\maketitle

\begin{quote}
{\bf Abstract}
In this paper, we study the construction methods for uninorms on bounded lattices via functions with the given uninorms and $q\in \mathbb{L_{B}}$ (or $p\in \mathbb{L_{B}}$). Specifically,  we investigate the conditions under which these functions can be uninorms on bounded lattices  when $q\in (0,\mathfrak{e})\cup I_{\mathfrak{e}}^{\varrho}$ and  $q\in  I_{\varrho}^{\mathfrak{e}}$    (or $p\in (\mathfrak{e},1)\cup I_{\mathfrak{e}}^{\sigma}$ and $p\in I_{\sigma}^{\mathfrak{e}}$), respectively. Moreover, some illustrative examples and figures are provided.

 %under some conditions on bounded lattices and the given uninorms with $q\in (0,e)\cup I_{e}^{a}\cup I_{a}^{e}$ and $q\in (e,1)\cup I_{e}^{b}\cup I_{b}^{e}$, respectively.
%In this paper, we study the construction methods for uninorms on bounded lattices via functions with uninorms and $q$ or $p$. And we illustrate that the functions given by ($\ref{eq1}$) and ($\ref{eq2}$) can be uninorms on bounded lattices under some conditions on bounded lattices and the given uninorms with $q\in (0,e)\cup I_{e}^{a}\cup I_{a}^{e}$ and $q\in (e,1)\cup I_{e}^{b}\cup I_{b}^{e}$, respectively. Moreover, some illustrative examples are provided.

{\textbf{Keywords}:}\ bounded lattices; $t$-norms; $t$-conorms; uninorms
\end{quote}

\section{Introduction}\label{intro}
The concept of triangular norms ($t$-norms for short) and triangular conorms ($t$-conorms  for short) on $[0,1]$, introduced by Menger \cite{KM42}, Schweizer and Sklar \cite{BS60,BS83}, have been widely used in several aspects, such as fuzzy logic, fuzzy set theory and so on (see, e.g., \cite{MB08,UH95,XL09,BS61,ZW07,HJ01}). And then, as the generalization of $t$-norms and $t$-conorms, the notion of uninorms on $[0,1]$, proposed by Yager and Rybalov \cite{RR96}, is an important tool in many fields, such as fuzzy logics, neural networks, expert systems and so on (see, e.g., \cite{BD99,MG09,MG11,WP07}).

Recently,  the construction methods for the above operators have been widely investigated on bounded lattices  by many researchers.
%, especially, the construction methods for them.
In fact, the researchers proposed  construction methods for $t$-norms (resp. $t$-conorms) via  $t$-norms (resp. $t$-conorms) (see, e.g., \cite{GD018,GD19.3,UE19,UE20,UE15,FK19,JM12,SS06,SS08}), $t$-subnorms (resp. $t$-superconorms) (see, e.g., \cite{ZX23}) and closure operators (resp. interior operators) (see, e.g., \cite{EA22.1,EA22.2}).
Afterward, the researchers also provided  construction methods for uninorms  via $t$-norms (resp. $t$-conorms) (see, e.g., \cite{EA21,SB14,GD19.2,GD19,GD19.1,GD20,GD17,GD16,YD19,YD20,FK15,AX20}), $t$-subnorms (resp. $t$-superconorms) (see, e.g., \cite{XJ20,WJ21,HP21}), closure operators (resp. interior operators) (see, e.g., \cite{GD21,XJ21,YO20,BZ21}), additive  generators (see, e.g., \cite{HeP}) and uninorms (see, e.g., \cite{GD23,ZX23,ZX24}).

More especially, in \cite{GD23} and \cite{ZX23}, the researchers  proposed new construction methods for uninorms on bounded lattices via a given uninorm on $[0,\varrho]$ (or $[\sigma,1]$) of $ \mathbb{L_{B}}$, respectively.  These  methods both stem  from  given uninorms on  $\mathbb{L_{B}}$
 % provided a novel perspective to study the constructions mthods for uninorms on bounded lattices
  and generalize some methods for uninorms, $t$-norms and $t$-conorms  in the literature.
In \cite{GD23}, G.D. \c{C}ayl{\i} et al. gave some methods for uninorms with the uninorms on $[0,\varrho]$ and  $t$-conorms on $[\varrho,1]$ (or uninorms on $[\sigma,1]$ and $t$-norms on $[0,\sigma]$) on bounded lattices under the condition that  $\varrho<\iota$ for $\iota \in \mathbb{L_{B}} \setminus [0,\varrho]$ (or the condition that $\iota <\sigma$ for $x\in \mathbb{L_{B}}\setminus [\sigma,1]$).
In \cite{ZX23}, Xiu and Zheng proposed the new construction method for uninorms on arbitrary bounded lattices with the given uninorm $\mathbb{U}^{*}$ satisfying the condition that  $\mathbb{U}^{*}(\iota, \kappa)\in[0,\mathfrak{e}]$ implies $(\iota , \kappa)\in [0,\mathfrak{e}]^{2}$ on $[0,\varrho]$ and the $t$-superconorms $\mathbb{R}$ with $\mathbb{R}(\varrho, \varrho)=\varrho$ on $[\varrho,1]$ (or  the given uninorm $\mathbb{U}^{*}$ satisfying  the condition that $\mathbb{U}^{*}(\iota, \kappa)\in[\mathfrak{e},1]$ implies $(\iota , \kappa)\in [\mathfrak{e},1]^{2}$ on $[\sigma,1]$) and the $t$-subnorms $\mathbb{F}$ with $\mathbb{F}(\sigma,\varrho)=\sigma$ on $[0,\sigma]$).
%In \cite{ZX23}, Xiu and Zheng proposed the construction methods for uninorms on arbitrary bounded lattices with the given uninorm $U^{*}$ satisfying the condition that  $U^{*}(x, y)\in[0,e]$ implies $(x,y)\in [0,e]^{2}$ on $[0,a]$ (or the condition that $U^{*}(x, y)\in[e,1]$ implies $(x,y)\in [e,1]^{2}$ on $[b,1]$) and the $t$-superconorms $R$ with $R(a,a)=a$ on $[a,1]$ (or the $t$-subnorms $F$ with $F(b,a)=b$ on $[0,b]$).
%In fact, the above methods both stem  from  given uninorms on  $L$
% % provided a novel perspective to study the constructions mthods for uninorms on bounded lattices
%  and generalize some methods for uninorms, $t$-norms and $t$-conorms  in the literature.

Recently, in \cite{ZX24},
Xiu and Zheng  proposed  the function with a given uninorm and $q\in \mathbb{L_{B}}$ (or $p\in \mathbb{L_{B}}$) and then  discussed how the function can be a uninorm with $q\in \{0\}\cup I_{\mathfrak{e}, \varrho}\cup [\mathfrak{e}, \varrho]\cup (\varrho,1]$ (or $p\in \{1\}\cup I_{\mathfrak{e},\sigma}\cup [\sigma,\mathfrak{e}]\cup [0,\sigma)$).
%also introduced some construction methods for uninorms on bounded lattices via uninorms under some conditions on $L$ and the given operators. Especially,
% Xiu and Zheng proposed two functions with   the given uninorms and $q$ (or $p$) and then
% $U$ in $(1)$ and $(2)$ on bounded lattices. And they
%discussed how the functions
%% given by $(\ref{eq1})$  and $(\ref{eq2})$)
%can be  uninorms with $q\in \{0\}\cup I_{e,a}\cup [e,a]\cup (a,1]$ and $p\in \{1\}\cup I_{e,b}\cup [b,e]\cup [0,b)$, respectively.
Motivated by this, we  study the functions with the given uninorm and $q$ (or $p$) on bounded lattices in this paper.
 %We first propose the function $U$ with a uninorm and a map with $p$ (or a map with $q$) in $(\ref{eq1})$ (or $(\ref{eq2})$) on bounded lattices. And then,
That is,  we try to find  conditions under which  the function  $(\ref{eq1})$ (or $(\ref{eq2})$)  with  the given uninorm and  $q\in (0,\mathfrak{e})\cup I_{\mathfrak{e}}^{\varrho}\cup I_{\varrho}^{\mathfrak{e}}$ (or $p\in (\mathfrak{e},1)\cup I_{\mathfrak{e}}^{\sigma}\cup I_{\sigma}^{\mathfrak{e}}$)
can be a uninorm on bounded lattices. In this case, besides the results in \cite{ZX24}, we have discussed all the cases of $q\in \mathbb{L_{B}}$ (or $p\in \mathbb{L_{B}}$) with which  the function  $(\ref{eq1})$ (or $(\ref{eq2})$)  can be a uninorm on bounded lattices.
In details, in Section 3, based on the given uninorm and $q\in \mathbb{L_{B}}$, we first define the functions $U$ by (\ref{eq1}).
% and (\ref{eq2}).
Then, we discuss how  the function  given by (\ref{eq1}) with $q\in (0,\mathfrak{e})\cup I_{\mathfrak{e}}^{\varrho}$ and  $q\in I_{\varrho}^{\mathfrak{e}}$ can be a uninorm, respectively.
Moreover, the dual results are also provided.  Meanwhile, we give some examples  and figures to illustrate the construction methods.

%the functions  given by (\ref{eq1}) and (\ref{eq2}) with all cases of $q$ and $p$ in $L$  have been discussed, respectively.

%given by (\ref{eq1}) (or (\ref{eq2})) with $q\in (0,e)\cup I_{e}^{a}\cup I_{a}^{e}$ (or $p\in (e,1)\cup I_{e}^{b}\cup I_{b}^{e}$) and then find some conditions under which the function $U$ is a uninorm on bounded lattices.

%The main results  of this paper is organized in Section 3.
%%In Section 2, we provide some basic notions and results on bounded lattices which will be used in this paper.
%In Section 3, based on the given uninorm and $q\in \mathbb{L_{B}}$, we first define the functions $U$ by (\ref{eq1}).
%% and (\ref{eq2}).
%Then, we discuss how  the function  given by (\ref{eq1}) with $q\in (0,\mathfrak{e})\cup I_{\mathfrak{e}}^{\varrho}$ and  $q\in I_{\varrho}^{\mathfrak{e}}$ can be a uninorm, respectively.
%Moreover, the dual results are also provided.  Meanwhile, we give some examples  and figures to illustrate the construction methods.
% In Section 4,  we end with some conclusions.

\section{Preliminaries}
In this section, we recall some basic notions of lattices and aggregation operators on bounded lattices.

\begin{definition}[\cite{GB67}]\label{de2.1}
A lattice $(\mathbb{L_{B}},\leq)$  is bounded if it has top  and bottom elements, which are written as $1$  and $0$, respectively, that is,  $0\leq \iota \leq 1$ for all $\iota \in \mathbb{L_{B}}$.
\end{definition}

%Throughout this article, unless stated otherwise, we denote $\mathbb{L_{B}}$  as a bounded lattice with the top and bottom elements  $1$ and $0$, respectively.

Throughout this article, $\mathbb{L_{B}}$ is denoted as a bounded lattice with the top element 1 and the  bottom element 0.

\begin{definition}[\cite{GB67}]\label{de2.2}
Let  $\varrho ,\sigma \in \mathbb{L_{B}}$ with $\varrho \leq \sigma$. A subinterval $[\varrho ,\sigma]$ of $\mathbb{L_{B}}$ is defined as
$$[\varrho , \sigma]=\{\iota \in \mathbb{L_{B}}: \varrho \leq \iota \leq \sigma\}.$$
Similarly, we can define $[\varrho , \sigma)=\{\iota \in \mathbb{L_{B}}: \varrho \leq \iota < \sigma\}, (\varrho , \sigma]=\{\iota\in \mathbb{L_{B}}: \varrho < \iota \leq \sigma\}$ and $(\varrho , \sigma)=\{\iota\in \mathbb{L_{B}}: \varrho < \iota < \sigma\}$. If $\varrho$  and $\sigma$  are incomparable, then we use the notation $\varrho \parallel \sigma$. If $\varrho$  and $\sigma$  are comparable, then we use the notation $\varrho \nparallel \sigma$.

\end{definition}

In the following,  $I_{\varrho}=\{\iota \in \mathbb{L_{B}}\mid \iota \parallel \varrho\}$.   $I^{\varrho}=\{\iota \in \mathbb{L_{B}} \mid \iota \nparallel \varrho \}$.
$I_{\varrho}^{\sigma}=\{\iota \in \mathbb{L_{B}}\mid \iota \parallel \varrho \ \mbox{and}\ \iota \nparallel \sigma\}$.  $I_{\varrho ,\sigma}=\{\iota \in \mathbb{L_{B}} \mid \iota \parallel \varrho \ \mbox{and} \ \iota \parallel \sigma\}$. Obviously, $I_{\varrho}^{\varrho}=\emptyset$ and $I_{\varrho,\varrho}=I_{\varrho}$.

\begin{definition}[\cite{SS06}]\label{de2.3}
 An operation $\mathbb{T}:\mathbb{L_{B}}^{2}\rightarrow \mathbb{L_{B}}$  is called a $t$-norm on $\mathbb{L_{B}}$ if it is commutative, associative, and increasing with respect to both variables, and it has the neutral element $1\in \mathbb{L_{B}}$.
\end{definition}

\begin{definition}[\cite{GD16}]\label{}
 An operation $\mathbb{S}:\mathbb{L_{B}}^{2}\rightarrow \mathbb{L_{B}}$  is called a $t$-conorm on  $\mathbb{L_{B}}$ if it is commutative, associative, and increasing with respect to both variables, and it has the neutral element $0\in \mathbb{L_{B}}$.
\end{definition}

\begin{definition}[\cite{FK15}]\label{de2.4}
 An operation $\mathbb{U}:\mathbb{L_{B}}^{2}\rightarrow \mathbb{L_{B}}$  is called a uninorm on  $\mathbb{L_{B}}$ (a uninorm if  $\mathbb{L_{B}}$ is fixed) if it is commutative, associative, and increasing with respect to both variables, and it has the neutral element $\mathfrak{e} \in \mathbb{L_{B}}$.
\end{definition}

\begin{proposition}[\cite{FK15}]\label{de34}
Let $\mathbb{U}$ be a uninorm on $\mathbb{L_{B}}$ with  $\mathfrak{e} \in \mathbb{L_{B}}\setminus\{0,1\}$. Then the following statements hold:\\
$(1)$ $\mathbb{T}_{e}=\mathbb{U}\mid [0,\mathfrak{e}]^{2}\rightarrow [0,\mathfrak{e}]$ is a $t$-norm on $[0,\mathfrak{e}]$.\\
$(2)$ $\mathbb{S}_{e}=\mathbb{U}\mid [\mathfrak{e},1]^{2}\rightarrow [\mathfrak{e},1]$ is a $t$-conorm on $[\mathfrak{e},1]$.
\end{proposition}

\begin{definition}[\cite{HP21}]\label{de35}
Let  $\mathfrak{e} \in \mathbb{L_{B}} \setminus\{0,1\}$. We denote by $\mathcal{U}_{min}$ the class of all uninorms $\mathbb{U}$ on $\mathbb{L_{B}}$ with neutral element $\mathfrak{e}$ satisfying $\mathbb{U}(\iota ,\kappa)=\kappa$, for all $(\iota ,\kappa)\in(\mathfrak{e},1]\times (\mathbb{L_{B}}\setminus[\mathfrak{e},1])$.

Similarly, we denote by $\mathcal{U}_{max}$ the class of all uninorms $\mathbb{U}$ on $\mathbb{L_{B}}$ with neutral element $\mathfrak{e}$ satisfying
$\mathbb{U}(\iota ,\kappa)=\kappa$, for all $(\iota ,\kappa)\in[0,\mathfrak{e})\times (\mathbb{L_{B}} \setminus[0,\mathfrak{e}])$.
\end{definition}

\begin{proposition}[\cite{WJ21}]\label{pro2.1}
Let $\mathcal{S}$ be a nonempty set and $\mathcal{C}_{1},\mathcal{C}_{2},\ldots,\mathcal{C}_{n}$ be subsets of $\mathcal{S}$. Let $\mathbb{G}$ be a commutative binary operation on $\mathcal{S}$. Then $\mathbb{G}$ is associative on $\mathcal{C}_{1}\cup \mathcal{C}_{2}\cup\ldots\cup \mathcal{C}_{n}$ if and only if all of the following statements hold:\\
$(i)$ for every combination $\{i,j,k\}$ of size $3$ chosen from $\{1,2,\ldots,n\}$, $\mathbb{G}(\iota,\mathbb{G}(\kappa , \omega ))=\mathbb{G}(\mathbb{G}(\iota,\kappa), \omega )=\mathbb{G}(\kappa,\mathbb{G}(\iota,\omega))$ for all $\iota \in \mathcal{C}_{i},\kappa \in \mathcal{C}_{j}, \omega \in \mathcal{C}_{k}$;\\
$(ii)$ for every combination $\{i,j\}$ of size $2$ chosen from $\{1,2,\ldots,n\}$, $\mathbb{G}(\iota,\mathbb{G}(\kappa ,\omega ))=\mathbb{G}(\mathbb{G}(\iota, \kappa ), \omega ) $ for all $\iota \in \mathcal{C}_{i},\kappa \in \mathcal{C}_{i}, \omega \in \mathcal{C}_{j}$;\\
$(iii)$ for every combination $\{i,j\}$ of size $2$ chosen from $\{1,2,\ldots,n\}$, $\mathbb{G}(\iota,\mathbb{G}(\kappa , \omega ))=\mathbb{G}(\mathbb{G}(\iota, \kappa ), \omega ) $ for all $\iota\in \mathcal{C}_{i},\kappa \in \mathcal{C}_{j},\omega \in \mathcal{C}_{j}$;\\
$(iv)$ for every $i\in \{1,2,\ldots,n\}$, $\mathbb{G}(\iota, \mathbb{G}(\kappa ,\omega ))=\mathbb{G}(\mathbb{G}(\iota,\kappa ),\omega ) $ for all $\iota,\kappa ,\omega \in \mathcal{C}_{i}$.
\end{proposition}

\begin{theorem}[\cite{GD018}]\label{th021}
 Let $\varrho \in \mathbb{L_{B}} \setminus \{0,1\}$. If $\mathbb{V}$ is a $t$-norm on $[\varrho ,1]$ and $\mathbb{W}$ is a $t$-conorm on $[0,\varrho]$, then $\mathbb{T}:\mathbb{L_{B}}^{2}\rightarrow \mathbb{L_{B}}$  is a $t$-norm   and $\mathbb{S}:\mathbb{L_{B}}^{2}\rightarrow \mathbb{L_{B}}$ is a $t$-conorm on $\mathbb{L_{B}}$, where

 $\mathbb{T}(\iota,\kappa )=\begin{cases}
 \mathbb{V}(\iota,\kappa ) &\mbox{if } (\iota,\kappa )\in [\varrho, 1]^{2},\\
 \iota \wedge \kappa &\mbox{if } 1\in \{\iota ,\kappa \},\\
 0 &\mbox{} otherwise,

\end{cases}$\\
and

$\mathbb{S}(\iota ,\kappa )=\begin{cases}
 \mathbb{W}(\iota ,\kappa ) &\mbox{if } (\iota ,\kappa )\in [0, \varrho]^{2},\\
 \iota \vee \kappa &\mbox{if } 0\in \{\iota ,\kappa \},\\
 1 &\mbox{} otherwise.

\end{cases}$\

\end{theorem}

\section{New construction methods for uninorms on bounded lattices}

In this section,  we study  the construction methods for uninorms on  a bounded lattice $\mathbb{L_{B}}$ via functions with the given uninorms and $q$ (or $p$).   Specifically,  based on functions via a given uninorm $\mathbb{U}^{*}$ on the subinterval $[0,\varrho]$ (or $[\sigma,1]$) of $\mathbb{L_{B}}$ and $q$ (or $p$),
    we propose  new methods to construct uninorms on  $\mathbb{L_{B}}$     under some conditions on $\mathbb{L_{B}}$  and $\mathbb{U}^{*}$.

%In this paper, we study the construction methods for uninorms on bounded lattices via functions with the given uninorms and $q$ (or $p$). Specifically,  we investigate the conditions under which these functions can be uninorms on bounded lattices  when $q\in (0,e)\cup I_{e}^{a}$ and  $q\in  I_{a}^{e}$    (or $p\in (e,1)\cup I_{e}^{b}$ and $p\in I_{b}^{e}$), respectively.

%In this section, based on a given uninorm $\mathbb{U}^{*}$ on the subinterval $[0,a]$ (or $[b,1]$) of $L$, we propose some new methods to construct uninorms on a bounded lattice $L$  with some conditions on $L$  and $\mathbb{U}^{*}$.
  For convenience, we denote by
$\mathcal{U}_{b}$  the class of all
uninorms  $\mathbb{U}$ on $\mathbb{L_{B}}$ with neutral element $\mathfrak{e}$ satisfying $\mathbb{U}(\iota , \kappa )\in[0,\mathfrak{e}]  $ implies $(\iota ,\kappa )\in [0,\mathfrak{e}]^{2}$.  Similarly, we denote by $\mathcal{U}_{t}$ the class of all uninorms $\mathbb{U}$ on $\mathbb{L_{B}}$ with neutral element $\mathfrak{e}$ satisfying  $\mathbb{U}(\iota , \kappa )\in[\mathfrak{e},1]  $ implies $(\iota ,\kappa )\in [\mathfrak{e},1]^{2}$.

Let $\varrho \in \mathbb{L_{B}}\setminus\{0,1\}$, $q\in \mathbb{L_{B}}$ and
$\mathbb{U}^{*}$ be a uninorm on $[0, \varrho]$ with a neutral element $\mathfrak{e}$. We can define a function $U:\mathbb{L_{B}}^{2}\rightarrow \mathbb{L_{B}}$ by
\begin{flalign}\label{eq1}
\ \ \ \ &\ \mathbb{U}(\iota ,\kappa )=\begin{cases}
\mathbb{U}^{*}(\iota , \kappa ) &\mbox{if } (\iota ,\kappa )\in [0, \varrho]^{2},\\
\iota  &\mbox{if } (\iota ,\kappa )\in (\mathbb{L_{B}} \setminus[0, \varrho ])\times [0,\mathfrak{e}],\\
\kappa  &\mbox{if } (\iota ,\kappa )\in [0,\mathfrak{e}]\times (\mathbb{L_{B}} \setminus[0, \varrho]),\\
\iota \vee \kappa \vee q &\mbox{if } (\iota ,\kappa )\in I_{\mathfrak{e}, \varrho}\times I_{\mathfrak{e}, \varrho},\\
1 &\mbox{}otherwise.
\end{cases}&
\end{flalign}

\begin{remark}
The structure of the function $\mathbb{U}$ given by  (\ref{eq1}) is illustrated in Fig.1.
\end{remark}

\begin{minipage}{11pc}
\setlength{\unitlength}{0.75pt}\begin{picture}(600,220)
\put(30,36){\makebox(0,0)[l]{\footnotesize$0$}}
\put(116,29){\makebox(0,0)[l]{\footnotesize$\mathfrak{e}$}}
\put(191,29){\makebox(0,0)[l]{\footnotesize$\varrho$}}
\put(266,29){\makebox(0,0)[l]{\footnotesize$1$}}
\put(300,29){\makebox(0,0)[l]{\footnotesize$I_{\mathfrak{e}}^{\varrho}$}}
\put(375,29){\makebox(0,0)[l]{\footnotesize$I_{\varrho}^{\mathfrak{e}}$}}
\put(455,29){\makebox(0,0)[l]{\footnotesize$I_{\mathfrak{e},\varrho}$}}

\put(30,69){\makebox(0,0)[l]{\footnotesize$\mathfrak{e}$}}
\put(30,99){\makebox(0,0)[l]{\footnotesize$\varrho$}}
\put(30,129){\makebox(0,0)[l]{\footnotesize$1$}}
\put(25,149){\makebox(0,0)[l]{\footnotesize$I_{\mathfrak{e}}^{\varrho}$}}
\put(25,179){\makebox(0,0)[l]{\footnotesize$I_{\varrho}^{\mathfrak{e}}$}}
\put(20,209){\makebox(0,0)[l]{\footnotesize$I_{\mathfrak{e}, \varrho}$}}

\put(56,53){\makebox(0,0)[l]{\footnotesize$\mathbb{U}^{*}(\iota ,\kappa )$}}
\put(131,53){\makebox(0,0)[l]{\footnotesize$\mathbb{U}^{*}(\iota ,\kappa )$}}
\put(56,85){\makebox(0,0)[l]{\footnotesize$\mathbb{U}^{*}(\iota ,\kappa )$}}
\put(131,85){\makebox(0,0)[l]{\footnotesize$\mathbb{U}^{*}(\iota ,\kappa )$}}

\put(76,115){\makebox(0,0)[l]{\footnotesize$\kappa $}}
\put(56,145){\makebox(0,0)[l]{\footnotesize$\mathbb{U}^{*}(\iota ,\kappa )$}}
\put(76,175){\makebox(0,0)[l]{\footnotesize$\kappa $}}
\put(76,205){\makebox(0,0)[l]{\footnotesize$\kappa $}}

\put(150,115){\makebox(0,0)[l]{\footnotesize$1$}}
\put(131,145){\makebox(0,0)[l]{\footnotesize$\mathbb{U}^{*}(\iota ,\kappa )$}}
\put(150,175){\makebox(0,0)[l]{\footnotesize$1$}}
\put(150,205){\makebox(0,0)[l]{\footnotesize$1$}}

\put(225,55){\makebox(0,0)[l]{\footnotesize$\iota $}}
\put(225,85){\makebox(0,0)[l]{\footnotesize$1$}}
\put(225,115){\makebox(0,0)[l]{\footnotesize$1$}}
\put(225,145){\makebox(0,0)[l]{\footnotesize$1$}}
\put(225,175){\makebox(0,0)[l]{\footnotesize$1$}}
\put(225,205){\makebox(0,0)[l]{\footnotesize$1$}}

\put(281,55){\makebox(0,0)[l]{\footnotesize$\mathbb{U}^{*}(\iota ,\kappa )$}}
\put(281,85){\makebox(0,0)[l]{\footnotesize$\mathbb{U}^{*}(\iota ,\kappa )$}}
\put(305,115){\makebox(0,0)[l]{\footnotesize$1$}}
\put(281,145){\makebox(0,0)[l]{\footnotesize$\mathbb{U}^{*}(\iota ,\kappa )$}}
\put(305,175){\makebox(0,0)[l]{\footnotesize$1$}}
\put(305,205){\makebox(0,0)[l]{\footnotesize$1$}}

\put(375,55){\makebox(0,0)[l]{\footnotesize$\iota $}}
\put(375,85){\makebox(0,0)[l]{\footnotesize$1$}}
\put(375,115){\makebox(0,0)[l]{\footnotesize$1$}}
\put(375,145){\makebox(0,0)[l]{\footnotesize$1$}}
\put(375,175){\makebox(0,0)[l]{\footnotesize$1$}}
\put(375,205){\makebox(0,0)[l]{\footnotesize$1$}}

\put(450,55){\makebox(0,0)[l]{\footnotesize$\iota $}}
\put(450,85){\makebox(0,0)[l]{\footnotesize$1$}}
\put(450,115){\makebox(0,0)[l]{\footnotesize$1$}}
\put(450,145){\makebox(0,0)[l]{\footnotesize$1$}}
\put(450,175){\makebox(0,0)[l]{\footnotesize$1$}}
\put(431,205){\makebox(0,0)[l]{\footnotesize$\iota \vee \kappa \vee q$}}

\put(45,39){\line(0,1){180}}
\put(120,39){\line(0,1){180}}
\put(195,39){\line(0,1){180}}
\put(270,39){\line(0,1){180}}
\put(345,39){\line(0,1){180}}
\put(420,39){\line(0,1){180}}
\put(495,39){\line(0,1){180}}

\put(270,39){\line(1,0){225}}
\put(270,39){\line(-1,0){225}}

\put(270,71){\line(1,0){225}}
\put(270,71){\line(-1,0){225}}

\put(270,101){\line(1,0){225}}
\put(270,101){\line(-1,0){225}}

\put(270,131){\line(1,0){225}}
\put(270,131){\line(-1,0){225}}

\put(270,161){\line(1,0){225}}
\put(270,161){\line(-1,0){225}}

\put(270,191){\line(1,0){225}}
\put(270,191){\line(-1,0){225}}

\put(270,219){\line(1,0){225}}
\put(270,219){\line(-1,0){225}}

\put(150,0){\emph{Fig.1. The function $\mathbb{U}$ given by (\ref{eq1})}.}
\end{picture}
\end{minipage}\\

In the following, we discuss the function $\mathbb{U}$ given by (\ref{eq1}) with $q\in(0,\mathfrak{e})\cup I_{\mathfrak{e}}^{\varrho}$  and $q\in I_{\varrho}^{\mathfrak{e}}$, respectively.

First, we illustrate that the function $\mathbb{U}$ given by (\ref{eq1}) with $q\in(0,\mathfrak{e})\cup I_{\mathfrak{e}}^{\varrho}$ can be a uninorm under some conditions.

\begin{theorem}\label{th31}
Let $\varrho \in \mathbb{L_{B}} \setminus\{0,1\}$, $q\in(0,\mathfrak{e})\cup I_{\mathfrak{e}}^{\varrho}$,
$\mathbb{U}^{*}$ be a uninorm on $[0,\varrho]$ with a neutral element $\mathfrak{e}$ and $\mathbb{U}_{[0,\varrho]}^{1}$ be a function given by $(\ref{eq1})$. Suppose that $\iota \vee \kappa =1$ for all $\iota ,\kappa \in I_{\mathfrak{e},\varrho}$ with $\iota \neq \kappa $, and $\iota \vee q =1$ for all $\iota \in I_{\mathfrak{e},\varrho}\cap I_{q}$.

$(1)$ Let us assume that $\mathbb{U}^{*}\in \mathcal{U}_{b}$. Then $\mathbb{U}_{[0,\varrho]}^{1}$ is a uninorm on $\mathbb{L_{B}} $ with the neutral element $\mathfrak{e} \in \mathbb{L_{B}} $ if and only if $\iota \parallel \kappa$ for all $\iota \in I_{\mathfrak{e},\varrho}\cap I^{q}$ and $\kappa \in I_{\mathfrak{e}}^{\varrho} $ .

$(2)$ Moreover, let us assume that $I_{\varrho}^{\mathfrak{e}}\cup I_{\mathfrak{e},\varrho}\cup (\varrho,1)\neq\emptyset$. Then $\mathbb{U}_{[0,\varrho]}^{1}$ is a uninorm on $\mathbb{L_{B}}$ with the neutral element $\mathfrak{e} \in \mathbb{L_{B}} $ if and only if $\mathbb{U}^{*}\in \mathcal{U}_{b}$ and $\iota \parallel \kappa $ for all $\iota \in I_{\mathfrak{e},\varrho}\cap I^{q}$ and $\kappa \in I_{\mathfrak{e}}^{\varrho} $.

\end{theorem}

\begin{proof}
(1) Necessity. Let $\mathbb{U}_{[0,\varrho]}^{1}$ be a uninorm on $\mathbb{L_{B}}$ with a neutral element $\mathfrak{e}$. We prove that $\iota \parallel \kappa $ for all $\iota \in I_{\mathfrak{e},\varrho} \cap I^{q}$ and $\kappa \in I_{\mathfrak{e}}^{\varrho} $.

Assume that there exist $\iota \in I_{\mathfrak{e},\varrho} \cap I^{q}$ and $y\in I_{\mathfrak{e}}^{\varrho} $  such that $\iota \nparallel \kappa $, i.e., $\kappa <\iota $ and $q<\iota $. Then $\mathbb{U}_{[0,\varrho]}^{1}(\iota ,\kappa )=1$ and $\mathbb{U}_{[0,\varrho]}^{1}(\iota ,\iota )=\iota \vee q=\iota $. Since $\iota <1$, the  increasingness property of $\mathbb{U}_{[0,\varrho]}^{1}$ is violated. Thus $\iota \parallel \kappa $   for all $\iota \in I_{\mathfrak{e},\varrho}\cap I^{q} $ and $\kappa \in I_{\mathfrak{e}}^{\varrho} $.

Sufficiency. It is obvious that  $\mathbb{U}_{[0,\varrho]}^{1}$ is commutative  and  $\mathfrak{e}$  is the neutral element of $\mathbb{U}_{[0,\varrho]}^{1}$. Thus, we just  prove the increasingness and the associativity of $\mathbb{U}_{[0,\varrho]}^{1}$.

I. Increasingness: We prove that if $\iota \leq \kappa $, then $\mathbb{U}_{[0,\varrho]}^{1}(\iota , \omega )\leq \mathbb{U}_{[0,\varrho]}^{1}(\kappa , \omega )$ for all $\omega \in \mathbb{L_{B}}$. It is obvious that $\mathbb{U}_{[0,\varrho]}^{1}(\iota , \omega )\leq \mathbb{U}_{[0,\varrho]}^{1}(\kappa , \omega )$ if both $\iota $ and $\kappa $ belong to one of the  intervals $ [0,\mathfrak{e}], I_{\mathfrak{e}}^{\varrho}, (\mathfrak{e},\varrho], I_{\varrho}^{\mathfrak{e}}, I_{\mathfrak{e},\varrho}$ or $(\varrho,1]$ for all $\omega \in \mathbb{L_{B}}$. The residual proof can be split into all possible cases:

1. $\iota \in [0,\mathfrak{e}]$

\ \ \ 1.1. $\kappa \in I_{\mathfrak{e}}^{\varrho}\cup (\mathfrak{e},\varrho]$

\ \ \ \ \ \ 1.1.1. $\omega \in [0,\mathfrak{e}]\cup I_{\mathfrak{e}}^{\varrho}\cup (\mathfrak{e},\varrho]$

\ \ \ \ \ \ \ \ \ \ \ \ $\mathbb{U}_{[0,\varrho]}^{1}(\iota ,\omega )=\mathbb{U}^{*}(\iota , \omega )\leq \mathbb{U}^{*}(\kappa ,\omega )=\mathbb{U}_{[0,\varrho]}^{1}(\kappa , \omega )$

\ \ \ \ \ \ 1.1.2. $\omega \in I_{\varrho}^{\mathfrak{e}}\cup I_{\mathfrak{e},\varrho}\cup(\varrho,1]$

\ \ \ \ \ \ \ \ \ \ \ \ $\mathbb{U}_{[0,\varrho]}^{1}(\iota ,\omega )=\omega \leq 1=\mathbb{U}_{[0,\varrho]}^{1}(\kappa , \omega )$

\ \ \ 1.2. $\kappa \in I_{\varrho}^{\mathfrak{e}}\cup (\varrho,1]$

\ \ \ \ \ \ 1.2.1. $\omega \in [0,\mathfrak{e}]$

\ \ \ \ \ \ \ \ \ \ \ \ $\mathbb{U}_{[0,\varrho]}^{1}(\iota ,\omega )=\mathbb{U}^{*}(\iota ,\omega )\leq \iota <\kappa =\mathbb{U}_{[0,\varrho]}^{1}(\kappa ,\omega )$

\ \ \ \ \ \ 1.2.2. $\omega \in I_{\mathfrak{e}}^{\varrho}\cup (\mathfrak{e},\varrho]$

\ \ \ \ \ \ \ \ \ \ \ \ $\mathbb{U}_{[0,\varrho]}^{1}(\iota , \omega )=\mathbb{U}^{*}(\iota , \omega )\leq \varrho<1=\mathbb{U}_{[0,\varrho]}^{1}(\kappa ,\omega )$

\ \ \ \ \ \ 1.2.3. $\omega \in I_{\varrho}^{\mathfrak{e}}\cup I_{\mathfrak{e}, \varrho}\cup(\varrho,1]$

\ \ \ \ \ \ \ \ \ \ \ \ $\mathbb{U}_{[0,\varrho]}^{1}(\iota ,\omega )=\omega \leq 1=\mathbb{U}_{[0,\varrho]}^{1}(\kappa ,\omega)$

\ \ \ 1.3. $\kappa \in I_{\mathfrak{e},\varrho}$

\ \ \ \ \ \ 1.3.1. $\omega \in [0,\mathfrak{e}]$

\ \ \ \ \ \ \ \ \ \ \ \ $\mathbb{U}_{[0,\varrho]}^{1}(\iota ,\omega )=\mathbb{U}^{*}(\iota ,\omega )\leq \iota <\kappa =\mathbb{U}_{[0,\varrho]}^{1}(\kappa ,\omega )$

\ \ \ \ \ \ 1.3.2. $\omega \in I_{\mathfrak{e}}^{\varrho}\cup (\mathfrak{e},\varrho]$

\ \ \ \ \ \ \ \ \ \ \ \ $\mathbb{U}_{[0,\varrho]}^{1}(\iota ,\omega )=\mathbb{U}^{*}(\iota , \omega )\leq \varrho <1=\mathbb{U}_{[0,\varrho]}^{1}(\kappa ,\omega )$

\ \ \ \ \ \ 1.3.3. $\omega \in I_{\varrho}^{\mathfrak{e}} \cup(\varrho,1]$

\ \ \ \ \ \ \ \ \ \ \ \ $\mathbb{U}_{[0,\varrho]}^{1}(\iota ,\omega )=\omega \leq 1=\mathbb{U}_{[0,\varrho]}^{1}(\kappa ,\omega )$

\ \ \ \ \ \ 1.3.4. $\omega \in I_{\mathfrak{e},\varrho} $

\ \ \ \ \ \ \ \ \ \ \ \ $\mathbb{U}_{[0,\varrho]}^{1}(\iota ,\omega )=\omega \leq \kappa \vee \omega \vee q=\mathbb{U}_{[0,\varrho]}^{1}(\kappa ,\omega )$

2. $\iota \in I_{\mathfrak{e}}^{\varrho}$

\ \ \ 2.1. $\kappa \in (\mathfrak{e},\varrho]$

\ \ \ \ \ \ 2.1.1. $\omega \in [0,\mathfrak{e}]\cup I_{\mathfrak{e}}^{\varrho}\cup (\mathfrak{e},\varrho]$

\ \ \ \ \ \ \ \ \ \ \ \ $\mathbb{U}_{[0,\varrho]}^{1}(\iota ,\omega )=\mathbb{U}^{*}(\iota ,\omega )\leq \mathbb{U}^{*}(\kappa ,\omega )=\mathbb{U}_{[0,\varrho]}^{1}(\kappa ,\omega )$

\ \ \ \ \ \ 2.1.2. $\omega \in I_{\varrho}^{\mathfrak{e}}\cup I_{\mathfrak{e}, \varrho}\cup (\varrho, 1]$

\ \ \ \ \ \ \ \ \ \ \ \ $\mathbb{U}_{[0,\varrho]}^{1}(\iota ,\omega )=1=\mathbb{U}_{[0,\varrho]}^{1}(\kappa ,\omega )$

\ \ \ 2.2. $\kappa \in I_{\varrho}^{\mathfrak{e}} \cup(\varrho, 1]$

\ \ \ \ \ \ 2.2.1. $\omega \in [0,\mathfrak{e}]$

\ \ \ \ \ \ \ \ \ \ \ \ $\mathbb{U}_{[0,\varrho]}^{1}(\iota ,\omega )=\mathbb{U}^{*}(\iota ,\omega )\leq \iota < \kappa =\mathbb{U}_{[0,\varrho]}^{1}(\kappa ,\omega )$

\ \ \ \ \ \ 2.2.2. $\omega \in I_{\mathfrak{e}}^{\varrho}\cup (\mathfrak{e}, \varrho]$

\ \ \ \ \ \ \ \ \ \ \ \ $\mathbb{U}_{[0,\varrho]}^{1}(\iota ,\omega )=\mathbb{U}^{*}(\iota ,\omega )\leq \varrho <1=\mathbb{U}_{[0,\varrho]}^{1}(\kappa , \omega )$

\ \ \ \ \ \ 2.2.3. $\omega \in I_{\varrho }^{\mathfrak{e}}\cup I_{\mathfrak{e}, \varrho}\cup (\varrho,1]$

\ \ \ \ \ \ \ \ \ \ \ \ $\mathbb{U}_{[0,\varrho]}^{1}(\iota ,\omega )=1=\mathbb{U}_{[0,\varrho]}^{1}(\kappa ,\omega )$

\ \ \ 2.3. $\kappa \in I_{\mathfrak{e},\varrho}$

\ \ \ \ \ \ 2.3.1. $\omega \in [0,\mathfrak{e}]$

\ \ \ \ \ \ \ \ \ \ \ \ $\mathbb{U}_{[0,\varrho]}^{1}(\iota ,\omega )=\mathbb{U}^{*}(\iota ,\omega )\leq \iota < \kappa =\mathbb{U}_{[0,\varrho]}^{1}(\kappa ,\omega )$

\ \ \ \ \ \ 2.3.2. $\omega \in I_{\mathfrak{e}}^{\varrho}\cup (\mathfrak{e},\varrho]$

\ \ \ \ \ \ \ \ \ \ \ \ $\mathbb{U}_{[0,\varrho]}^{1}(\iota ,\omega )=\mathbb{U}^{*}(\iota , \omega )\leq \varrho<1=\mathbb{U}_{[0,\varrho]}^{1}(\kappa , \omega )$

\ \ \ \ \ \ 2.3.3. $\omega \in I_{\varrho}^{\mathfrak{e}}\cup (\varrho,1]$

\ \ \ \ \ \ \ \ \ \ \ \ $\mathbb{U}_{[0,\varrho]}^{1}(\iota ,\omega )=1=\mathbb{U}_{[0,\varrho]}^{1}(\kappa ,\omega )$

\ \ \ \ \ \ 2.3.4. $\omega \in I_{\mathfrak{e},\varrho}$

\ \ \ \ \ \ \ \ \ \ \ \ $\mathbb{U}_{[0,\varrho]}^{1}(\iota ,\omega )=1=\kappa \vee \omega \vee q = \mathbb{U}_{[0,\varrho]}^{1}(\kappa ,\omega )$

3. $\iota \in (\mathfrak{e},\varrho], \kappa \in I_{\varrho}^{\mathfrak{e}}\cup(\varrho,1]$

\ \ \ 3.1. $\omega \in [0,\mathfrak{e}]$

\ \ \ \ \ \ \ \ \ \ \ \ $\mathbb{U}_{[0,\varrho]}^{1}(\iota ,\omega )=\mathbb{U}^{*}(\iota ,\omega )\leq \iota <\kappa =\mathbb{U}_{[0,\varrho]}^{1}(\kappa ,\omega )$

\ \ \ 3.2. $\omega \in I_{\mathfrak{e}}^{\varrho}\cup (\mathfrak{e},\varrho]$

\ \ \ \ \ \ \ \ \ \ \ \ $\mathbb{U}_{[0,\varrho]}^{1}(\iota ,\omega )=\mathbb{U}^{*}(\iota ,\omega )\leq \varrho <1=\mathbb{U}_{[0,\varrho]}^{1}(\kappa ,\omega )$

\ \ \ 3.3. $\omega \in I_{\varrho}^{\mathfrak{e}}\cup I_{\mathfrak{e},\varrho}\cup(\varrho,1]$

\ \ \ \ \ \ \ \ \ \ \ \ $\mathbb{U}_{[0,\varrho]}^{1}(\iota ,\omega )=1=\mathbb{U}_{[0,\varrho]}^{1}(\kappa ,\omega )$

4. $\iota \in I_{\varrho}^{\mathfrak{e}}, \kappa \in (\varrho,1]$

\ \ \ 4.1. $\omega \in [0,\mathfrak{e}]$

\ \ \ \ \ \ \ \ \ \ \ \ $\mathbb{U}_{[0,\varrho]}^{1}(\iota ,\omega )=\iota \leq \kappa =\mathbb{U}_{[0,\varrho]}^{1}(\kappa ,\omega )$

\ \ \ 4.2. $\omega \in I_{\mathfrak{e}}^{\varrho}\cup(\mathfrak{e}, \varrho]\cup I_{\varrho}^{\mathfrak{e}}\cup I_{\mathfrak{e}, \varrho}\cup(\varrho,1]$

\ \ \ \ \ \ \ \ \ \ \ \ $\mathbb{U}_{[0,\varrho]}^{1}(\iota ,\omega )=1=\mathbb{U}_{[0,\varrho]}^{1}(\kappa ,\omega )$

5. $\iota \in I_{\mathfrak{e},\varrho}, \kappa \in I_{\varrho}^{\mathfrak{e}}\cup(\varrho,1]$

\ \ \ 5.1. $\omega \in [0,\mathfrak{e}]$

\ \ \ \ \ \ \ \ \ \ \ \ $\mathbb{U}_{[0,\varrho]}^{1}(\iota ,\omega )=\iota < \kappa =\mathbb{U}_{[0,\varrho]}^{1}(\kappa ,\omega )$

\ \ \ 5.2. $\omega \in I_{\mathfrak{e}}^{\varrho}\cup(\mathfrak{e},\varrho]\cup I_{\varrho}^{\mathfrak{e}}\cup (\varrho,1]$

\ \ \ \ \ \ \ \ \ \ \ \ $\mathbb{U}_{[0,\varrho]}^{1}(\iota ,\omega )=1=\mathbb{U}_{[0,\varrho]}^{1}(\kappa ,\omega )$

\ \ \ 5.3. $\omega \in I_{\mathfrak{e}, \varrho} $

\ \ \ \ \ \ \ \ \ \ \ \ $\mathbb{U}_{[0,\varrho]}^{1}(\iota ,\omega )=\iota \vee \omega \vee q\leq 1=\mathbb{U}_{[0,\varrho]}^{1}(\kappa ,\omega )$

II. Associativity: We demonstrate that $\mathbb{U}_{[0,\varrho]}^{1}(\iota ,\mathbb{U}_{[0,\varrho]}^{1}(\kappa ,\omega ))=\mathbb{U}_{[0,\varrho]}^{1}(U_{1}(\iota ,\kappa ),\omega )$ for all $\iota ,\kappa ,\omega \in \mathbb{L_{B}}$. By Proposition \ref{pro2.1}, we need to consider the following cases:

1. If $\iota ,\kappa ,\omega \in [0,\mathfrak{e}]\cup I_{\mathfrak{e}}^{\varrho}\cup (\mathfrak{e}, \varrho]$, then since $\mathbb{U}^{*}$ is associative, we have $\mathbb{U}_{[0,\varrho]}^{1}(\iota ,\mathbb{U}_{[0,\varrho]}^{1}(\kappa ,\omega ))=\mathbb{U}_{[0,\varrho]}^{1}(\mathbb{U}_{[0,\varrho]}^{1}(\iota ,\kappa ), \omega )=\mathbb{U}_{[0,\varrho]}^{1}(\kappa , \mathbb{U}_{[0,\varrho]}^{1}(\iota ,\omega ))$.

2. If $\iota ,\kappa ,\omega \in I_{\varrho}^{\mathfrak{e}}\cup (\varrho,1]$, then $\mathbb{U}_{[0,\varrho]}^{1}(\iota ,\mathbb{U}_{[0,\varrho]}^{1}(\kappa , \omega ))=\mathbb{U}_{[0,\varrho]}^{1}(\kappa ,\mathbb{U}_{[0,\varrho]}^{1}(\iota , \omega ))=\mathbb{U}_{[0,\varrho]}^{1}(\mathbb{U}_{[0,\varrho]}^{1}(\iota , \kappa ), \omega )=1$.

3. Assume that $\iota ,\kappa , \omega \in  I_{\mathfrak{e}, \varrho} $.

3.1. Suppose that $\iota ,\kappa , \omega \nparallel q$.

3.1.1. If $\iota \neq \kappa $, $\kappa \neq \omega $ and $\iota \neq \omega $, then $\mathbb{U}_{[0,\varrho]}^{1}(\iota ,\mathbb{U}_{[0,\varrho]}^{1}(\kappa , \omega ))=\mathbb{U}_{[0,\varrho]}^{1}(\iota ,\kappa \vee \omega \vee q)=\mathbb{U}_{[0,\varrho]}^{1}(\iota ,1)=1=\mathbb{U}_{[0,\varrho]}^{1}(\iota \vee \kappa \vee q,\omega )= \mathbb{U}_{[0,\varrho]}^{1}(\mathbb{U}_{[0,\varrho]}^{1}(\iota ,\kappa ),\omega )$  and $\mathbb{U}_{[0,\varrho]}^{1}(\kappa ,\mathbb{U}_{[0,\varrho]}^{1}(\iota ,\omega ))=\mathbb{U}_{[0,\varrho]}^{1}(\kappa ,\iota \vee \omega \vee q)=\mathbb{U}_{[0,\varrho]}^{1}(\kappa ,1)=1$.

3.1.2. If $\iota =\kappa $ and $\iota ,\kappa \neq \omega $, then $\mathbb{U}_{[0,\varrho]}^{1}(\iota ,\mathbb{U}_{[0,\varrho]}^{1}(\kappa ,\omega ))=\mathbb{U}_{[0,\varrho]}^{1}(\iota ,\kappa \vee \omega \vee q)=\mathbb{U}_{[0,\varrho]}^{1}(\iota ,1)=1=\iota \vee \omega \vee q=\mathbb{U}_{[0,\varrho]}^{1}(\iota ,\omega )= \mathbb{U}_{[0,\varrho]}^{1}(\mathbb{U}_{[0,\varrho]}^{1}(\iota ,\iota ),\omega )= \mathbb{U}_{[0,\varrho]}^{1}(\mathbb{U}_{[0,\varrho]}^{1}(\iota ,\kappa ),\omega )$  and $\mathbb{U}_{[0,\varrho]}^{1}(\kappa ,\mathbb{U}_{[0,\varrho]}^{1}(\iota ,\omega ))=\mathbb{U}_{[0,\varrho]}^{1}(\iota ,\mathbb{U}_{[0,\varrho]}^{1}(\iota ,\omega ))=\mathbb{U}_{[0,\varrho]}^{1}(\iota ,\mathbb{U}_{[0,\varrho]}^{1}(\kappa ,\omega ))=1$.

3.1.3. If  $\kappa = \omega $ and $\kappa ,\omega \neq \iota $, then we also have $\mathbb{U}_{[0,\varrho]}^{1}(\iota ,\mathbb{U}_{[0,\varrho]}^{1}(\kappa ,\omega ))=\mathbb{U}_{[0,\varrho]}^{1}(\mathbb{U}_{[0,\varrho]}^{1}(\iota ,\kappa ),\omega )=\mathbb{U}_{[0,\varrho]}^{1}(\kappa ,\mathbb{U}_{[0,\varrho]}^{1}(\iota ,\omega ))$ by the commutativity property of $\mathbb{U}_{[0,\varrho]}^{1}$.

3.1.4. If $\iota = \omega $ and $\iota ,\omega \neq \kappa$, then we also have $\mathbb{U}_{[0,\varrho]}^{1}(\iota ,\mathbb{U}_{[0,\varrho]}^{1}(\kappa ,\omega ))= \mathbb{U}_{[0,\varrho]}^{1}(\mathbb{U}_{[0,\varrho]}^{1}(\iota ,\kappa ),\omega )=\mathbb{U}_{[0,\varrho]}^{1}(\kappa ,\mathbb{U}_{[0,\varrho]}^{1}(\iota ,\omega ))$ by the commutativity property of $\mathbb{U}_{[0,\varrho]}^{1}$.

3.1.5. If $\iota = \kappa  = \omega $, then we can easily obtain $\mathbb{U}_{[0,\varrho]}^{1}(\iota ,\mathbb{U}_{[0,\varrho]}^{1}(\kappa ,\omega ))= \mathbb{U}_{[0,\varrho]}^{1}(\mathbb{U}_{[0,\varrho]}^{1}(\iota ,\kappa ),\omega )=\mathbb{U}_{[0,\varrho]}^{1}(\kappa ,\mathbb{U}_{[0,\varrho]}^{1}(\iota ,\omega ))$.

3.2. Suppose that there exist $ \iota \in I_{\mathfrak{e},\varrho}$ such that $\iota \parallel q$.

3.2.1. If $\iota \nparallel q$ and $\kappa ,\omega \parallel q$, then $\mathbb{U}_{[0,\varrho]}^{1}(\iota ,\mathbb{U}_{[0,\varrho]}^{1}(\kappa ,\omega ))=\mathbb{U}_{[0,\varrho]}^{1}(\iota ,\kappa \vee \omega  \vee q)=\mathbb{U}_{[0,\varrho]}^{1}(\iota ,1)=1=\mathbb{U}_{[0,\varrho]}^{1}(1,\omega )=\mathbb{U}_{[0,\varrho]}^{1}(\iota \vee \kappa \vee q,\omega )= \mathbb{U}_{[0,\varrho]}^{1}(\mathbb{U}_{[0,\varrho]}^{1}(\iota ,\kappa ),\omega )$ and $\mathbb{U}_{[0,\varrho]}^{1}(\kappa ,\mathbb{U}_{[0,\varrho]}^{1}(\iota ,\omega ))=\mathbb{U}_{[0,\varrho]}^{1}(\kappa ,\iota \vee \omega \vee q)=\mathbb{U}_{[0,\varrho]}^{1}(\kappa ,1)=1$.

3.2.2. If $\iota ,\kappa \nparallel q$ and $\omega \parallel q$, then $\mathbb{U}_{[0,\varrho]}^{1}(\iota ,\mathbb{U}_{[0,\varrho]}^{1}(\kappa ,\omega ))=\mathbb{U}_{[0,\varrho]}^{1}(\iota ,1)=1=\mathbb{U}_{[0,\varrho]}^{1}(\kappa ,\iota \vee \omega \vee q)= \mathbb{U}_{[0,\varrho]}^{1}(\kappa ,\mathbb{U}_{[0,\varrho]}^{1}(\iota , \omega ))$ and $\mathbb{U}_{[0,\varrho]}^{1}(\mathbb{U}_{[0,\varrho]}^{1}(\iota ,\kappa ),\omega )=\mathbb{U}_{[0,\varrho]}^{1}(\iota \vee \kappa \vee q,\omega )$. Moreover, we can obtain that if $\iota =\kappa $, then $\iota \vee \kappa \vee q=\iota $, $\mathbb{U}_{[0,\varrho]}^{1}(\mathbb{U}_{[0,\varrho]}^{1}(\iota ,\kappa ),\omega )=\mathbb{U}_{[0,\varrho]}^{1}(\iota ,\omega )=1$ and if $\iota \neq \kappa $, then $\iota \vee \kappa \vee q=1$, $\mathbb{U}_{[0,\varrho]}^{1}(\mathbb{U}_{[0,\varrho]}^{1}(\iota ,\kappa ),\omega )=\mathbb{U}_{[0,\varrho]}^{1}(1,\omega )=1$. Thus $\mathbb{U}_{[0,\varrho]}^{1}(\iota ,\mathbb{U}_{[0,\varrho]}^{1}(\kappa ,\omega )) =\mathbb{U}_{[0,\varrho]}^{1}(\mathbb{U}_{[0,\varrho]}^{1}(\iota ,\kappa ),\omega )=\mathbb{U}_{[0,\varrho]}^{1}(\kappa ,\mathbb{U}_{[0,\varrho]}^{1}(\iota ,\omega))$.

3.2.3. If $\iota , \kappa ,\omega \parallel q$, then $\mathbb{U}_{[0,\varrho]}^{1}(\iota ,\mathbb{U}_{[0,\varrho]}^{1}(\kappa ,\omega ))=\mathbb{U}_{[0,\varrho]}^{1}(\iota ,1)=1=\mathbb{U}_{[0,\varrho]}^{1}(1,\omega )=\mathbb{U}_{[0,\varrho]}^{1}(\iota \vee \kappa \vee q,\omega )= \mathbb{U}_{[0,\varrho]}^{1}(\mathbb{U}_{[0,\varrho]}^{1}(\iota ,\kappa ),\omega )$ and $\mathbb{U}_{[0,\varrho]}^{1}(\kappa ,\mathbb{U}_{[0,\varrho]}^{1}(\iota ,\omega ))=\mathbb{U}_{[0,\varrho]}^{1}(\kappa ,\iota \vee \kappa \vee q)=1$.

4. If $\iota ,\kappa \in [0,\mathfrak{e}]$ and $\omega \in I_{\varrho}^{\mathfrak{e}}\cup I_{\mathfrak{e},\varrho}\cup  (\varrho,1]$, then $\mathbb{U}_{[0,\varrho]}^{1}(\iota ,\mathbb{U}_{[0,\varrho]}^{1}(\kappa ,\omega ))=\mathbb{U}_{[0,\varrho]}^{1}(\iota ,\omega )=\omega =\mathbb{U}_{[0,\varrho]}^{1}(\mathbb{U}^{*}(\iota ,\kappa ),\omega)=\mathbb{U}_{[0,\varrho]}^{1}(\mathbb{U}_{[0,\varrho]}^{1}(\iota ,\kappa ),\omega) $.

5. If $\iota ,\kappa \in I_{\mathfrak{e}}^{\varrho}\cup (\mathfrak{e},\varrho]$ and $\omega \in I_{\varrho}^{\mathfrak{e}}\cup I_{\mathfrak{e},\varrho}\cup  (\varrho,1]$, then $\mathbb{U}_{[0,\varrho]}^{1}(\iota ,\mathbb{U}_{[0,\varrho]}^{1}(\kappa ,\omega ))=\mathbb{U}_{[0,\varrho]}^{1}(\iota ,1)=1=\mathbb{U}_{[0,\varrho]}^{1}(\mathbb{U}^{*}(\iota ,\kappa ),\omega )$ $=\mathbb{U}_{[0,\varrho]}^{1}(\mathbb{U}_{[0,\varrho]}^{1}(\iota ,\kappa ),\omega) $ and $\mathbb{U}_{[0,\varrho]}^{1}(\kappa ,\mathbb{U}_{[0,\varrho]}^{1}(\iota ,\omega ))=\mathbb{U}_{[0,\varrho]}^{1}(\kappa ,1)=1$. Thus $\mathbb{U}_{[0,\varrho]}^{1}(\iota ,\mathbb{U}_{[0,\varrho]}^{1}(\kappa ,\omega ))$ $=\mathbb{U}_{[0,\varrho]}^{1}(\mathbb{U}_{[0,\varrho]}^{1}(\iota ,\kappa ),\omega )=\mathbb{U}_{[0,\varrho]}^{1}(\kappa ,\mathbb{U}_{[0,\varrho]}^{1}(\iota ,\omega ))$.

6. If $\iota ,\kappa \in I_{\varrho}^{e}$ and $\omega \in  I_{\mathfrak{e},\varrho}\cup  (\varrho,1]$, then $\mathbb{U}_{[0,\varrho]}^{1}(\iota ,\mathbb{U}_{[0,\varrho]}^{1}(\kappa ,\omega ))=\mathbb{U}_{[0,\varrho]}^{1}(\iota ,1)=1=\mathbb{U}_{[0,\varrho]}^{1}(1,\omega )=\mathbb{U}_{[0,\varrho]}^{1}(\mathbb{U}_{[0,\varrho]}^{1}(\iota ,\kappa ),\omega )$.

7. If $\iota ,\kappa \in I_{\mathfrak{e},\varrho}$ and $\omega \in (\varrho,1]$, then $\mathbb{U}_{[0,\varrho]}^{1}(\iota ,\mathbb{U}_{[0,\varrho]}^{1}(\kappa ,\omega ))=\mathbb{U}_{[0,\varrho]}^{1}(\iota ,1)=1=\mathbb{U}_{[0,\varrho]}^{1}(\iota \vee \kappa \vee q,\omega )=\mathbb{U}_{[0,\varrho]}^{1}(\mathbb{U}_{[0,\varrho]}^{1}(\iota ,\kappa ),\omega )$.

8. If $\iota \in [0,\mathfrak{e}]$ and $\kappa ,\omega \in I_{\varrho}^{\mathfrak{e}}\cup (\varrho,1] $, then $\mathbb{U}_{[0,\varrho]}^{1}(\iota ,\mathbb{U}_{[0,\varrho]}^{1}(\kappa ,\omega ))=\mathbb{U}_{[0,\varrho]}^{1}(\iota ,1)=1=\mathbb{U}_{[0,\varrho]}^{1}(\kappa ,\omega )=\mathbb{U}_{[0,\varrho]}^{1}(\mathbb{U}_{[0,\varrho]}^{1}(\iota ,\kappa ),\omega )$ and $\mathbb{U}_{[0,\varrho]}^{1}(\kappa ,\mathbb{U}_{[0,\varrho]}^{1}(\iota ,\omega ))=\mathbb{U}_{[0,\varrho]}^{1}(\kappa ,\omega )=1$. Thus $\mathbb{U}_{[0,\varrho]}^{1}(\iota ,\mathbb{U}_{[0,\varrho]}^{1}(\kappa ,\omega ))=\mathbb{U}_{[0,\varrho]}^{1}(\mathbb{U}_{[0,\varrho]}^{1}(\iota ,\kappa ),\omega )=\mathbb{U}_{[0,\varrho]}^{1}(\kappa ,\mathbb{U}_{[0,\varrho]}^{1}(\iota ,\omega ))$.

9. If $\iota \in [0,\mathfrak{e}]$ and $\kappa ,\omega \in I_{\mathfrak{e},\varrho}$, then $\mathbb{U}_{[0,\varrho]}^{1}(\iota ,\mathbb{U}_{[0,\varrho]}^{1}(\kappa ,\omega ))=\mathbb{U}_{[0,\varrho]}^{1}(\iota ,\kappa \vee \omega \vee q)=\kappa \vee \omega \vee q=\mathbb{U}_{[0,\varrho]}^{1}(\kappa ,\omega )= \mathbb{U}_{[0,\varrho]}^{1}(\mathbb{U}_{[0,\varrho]}^{1}(\iota , \kappa ),\omega )$.

10. If $\iota \in I_{\mathfrak{e}}^{\varrho}\cup (\mathfrak{e},\varrho]$ and $\kappa ,\omega \in I_{\varrho}^{\mathfrak{e}} \cup (\varrho,1] $, then $\mathbb{U}_{[0,\varrho]}^{1}(\iota ,\mathbb{U}_{[0,\varrho]}^{1}(\kappa ,\omega ))=\mathbb{U}_{[0,\varrho]}^{1}(\iota ,1)=1=\mathbb{U}_{[0,\varrho]}^{1}(1,\omega )=\mathbb{U}_{[0,\varrho]}^{1}(\mathbb{U}_{[0,\varrho]}^{1}(\iota , \kappa ),\omega )$ and $\mathbb{U}_{[0,\varrho]}^{1}(\kappa ,\mathbb{U}_{[0,\varrho]}^{1}(\iota ,\omega ))=\mathbb{U}_{[0,\varrho]}^{1}(\kappa ,1)=1$. Thus $\mathbb{U}_{[0,\varrho]}^{1}(\iota ,\mathbb{U}_{[0,\varrho]}^{1}(\kappa ,\omega ))$ $=\mathbb{U}_{[0,\varrho]}^{1}(\mathbb{U}_{[0,\varrho]}^{1}(\iota , \kappa ),\omega )=\mathbb{U}_{[0,\varrho]}^{1}(\kappa , \mathbb{U}_{[0,\varrho]}^{1}(\iota ,\omega ))$.

11. If $\iota \in I_{\mathfrak{e}}^{\varrho}\cup (\mathfrak{e},\varrho]\cup I_{\varrho}^{\mathfrak{e}}$ and $\kappa ,\omega \in  I_{\mathfrak{e},\varrho}$, then $\mathbb{U}_{[0,\varrho]}^{1}(\iota ,\mathbb{U}_{[0,\varrho]}^{1}(\kappa ,\omega ))=\mathbb{U}_{[0,\varrho]}^{1}(\iota ,\kappa \vee \omega \vee q)=1=\mathbb{U}_{[0,\varrho]}^{1}(1,\omega )=\mathbb{U}_{[0,\varrho]}^{1}(\mathbb{U}_{[0,\varrho]}^{1}(\iota , \kappa ),\omega )$.

12. If $\iota \in I_{\mathfrak{e},\varrho}$ and $\kappa , \omega \in (\varrho,1]$, then $\mathbb{U}_{[0,\varrho]}^{1}(\iota ,\mathbb{U}_{[0,\varrho]}^{1}(\kappa ,\omega ))=\mathbb{U}_{[0,\varrho]}^{1}(\iota ,1)=1=\mathbb{U}_{[0,\varrho]}^{1}(1,\omega )=\mathbb{U}_{[0,\varrho]}^{1}(\mathbb{U}_{[0,\varrho]}^{1}(\iota , \kappa ),\omega )$.

13. If $\iota \in [0,\mathfrak{e}], \kappa \in I_{\mathfrak{e}}^{\varrho}\cup (\mathfrak{e},\varrho]$ and $\omega \in I_{\varrho}^{\mathfrak{e}}\cup I_{\mathfrak{e},\varrho}\cup  (\varrho,1] $, then $\mathbb{U}_{[0,\varrho]}^{1}(\iota ,\mathbb{U}_{[0,\varrho]}^{1}(\kappa ,\omega ))=\mathbb{U}_{[0,\varrho]}^{1}(\iota ,1)=1=\mathbb{U}_{[0,\varrho]}^{1}(\mathbb{U}^{*}(\iota , \kappa ),\omega )=\mathbb{U}_{[0,\varrho]}^{1}(\mathbb{U}_{[0,\varrho]}^{1}(\iota , \kappa ),\omega )$ and $\mathbb{U}_{[0,\varrho]}^{1}(\kappa ,\mathbb{U}_{[0,\varrho]}^{1}(\iota ,\omega ))=\mathbb{U}_{[0,\varrho]}^{1}(\kappa ,\omega )=1$. Thus $\mathbb{U}_{[0,\varrho]}^{1}(\iota ,\mathbb{U}_{[0,\varrho]}^{1}(\kappa ,\omega ))$ $=\mathbb{U}_{[0,\varrho]}^{1}(\mathbb{U}_{[0,\varrho]}^{1}(\iota , \kappa ),\omega )$ $=\mathbb{U}_{[0,\varrho]}^{1}(\kappa ,\mathbb{U}_{[0,\varrho]}^{1}(\iota ,\omega))$.

14. If $\iota \in [0,\mathfrak{e}],\kappa \in I_{\varrho}^{\mathfrak{e}}$ and $\omega \in I_{\mathfrak{e},\varrho}$, then $\mathbb{U}_{[0,\varrho]}^{1}(\iota ,\mathbb{U}_{[0,\varrho]}^{1}(\kappa ,\omega ))=\mathbb{U}_{[0,\varrho]}^{1}(\iota ,1)=1=\mathbb{U}_{[0,\varrho]}^{1}(\kappa ,\omega )=\mathbb{U}_{[0,\varrho]}^{1}(\mathbb{U}_{[0,\varrho]}^{1}(\iota ,\kappa ),\omega )$ and $\mathbb{U}_{[0,\varrho]}^{1}(\kappa ,\mathbb{U}_{[0,\varrho]}^{1}(\iota ,\omega ))=\mathbb{U}_{[0,\varrho]}^{1}(\kappa ,\omega )=1$. Thus $\mathbb{U}_{[0,\varrho]}^{1}(\iota ,\mathbb{U}_{[0,\varrho]}^{1}(\kappa ,\omega ))=\mathbb{U}_{[0,\varrho]}^{1}(\mathbb{U}_{[0,\varrho]}^{1}(\iota ,\kappa ),\omega)=\mathbb{U}_{[0,\varrho]}^{1}(\kappa ,\mathbb{U}_{[0,\varrho]}^{1}(\iota ,\omega))$.

15. If $\iota \in [0,\mathfrak{e}],\kappa \in I_{\mathfrak{e},\varrho}$ and $z\in (\varrho,1]$, then $\mathbb{U}_{[0,\varrho]}^{1}(\iota ,\mathbb{U}_{[0,\varrho]}^{1}(\kappa ,\omega))=\mathbb{U}_{[0,\varrho]}^{1}(\iota ,1)=1=\mathbb{U}_{[0,\varrho]}^{1}(\kappa ,\omega)=\mathbb{U}_{[0,\varrho]}^{1}(\mathbb{U}_{[0,\varrho]}^{1}(\iota ,\kappa ),\omega)$ and $\mathbb{U}_{[0,\varrho]}^{1}(\kappa ,\mathbb{U}_{[0,\varrho]}^{1}(\iota ,\omega))=\mathbb{U}_{[0,\varrho]}^{1}(\kappa ,\omega)=1$. Thus $\mathbb{U}_{[0,\varrho]}^{1}(\iota ,\mathbb{U}_{[0,\varrho]}^{1}(\kappa ,\omega))$ $=\mathbb{U}_{[0,\varrho]}^{1}(\mathbb{U}_{[0,\varrho]}^{1}(\iota ,\kappa ),\omega)=\mathbb{U}_{[0,\varrho]}^{1}(\kappa ,\mathbb{U}_{[0,\varrho]}^{1}(\iota ,\omega))$.

16. If $\iota \in I_{\mathfrak{e}}^{\varrho}\cup (\mathfrak{e},\varrho], \kappa \in I_{\varrho}^{\mathfrak{e}}$ and $z\in I_{\mathfrak{e},\varrho}$, then $\mathbb{U}_{[0,\varrho]}^{1}(\iota ,\mathbb{U}_{[0,\varrho]}^{1}(\kappa ,\omega))=\mathbb{U}_{[0,\varrho]}^{1}(\iota ,1)=1=\mathbb{U}_{[0,\varrho]}^{1}(1,z)=\mathbb{U}_{[0,\varrho]}^{1}(\mathbb{U}_{[0,\varrho]}^{1}(\iota , \kappa ),\omega )$ and $\mathbb{U}_{[0,\varrho]}^{1}(\kappa ,\mathbb{U}_{[0,\varrho]}^{1}(\iota ,z))=\mathbb{U}_{[0,\varrho]}^{1}(\kappa ,1)=1$. Thus $\mathbb{U}_{[0,\varrho]}^{1}(\iota ,\mathbb{U}_{[0,\varrho]}^{1}(\kappa ,\omega ))=\mathbb{U}_{[0,\varrho]}^{1}(\mathbb{U}_{[0,\varrho]}^{1}(\iota , \kappa ),\omega )=\mathbb{U}_{[0,\varrho]}^{1}(\kappa ,\mathbb{U}_{[0,\varrho]}^{1}(\iota ,\omega ))$.

17. If $\iota \in I_{\mathfrak{e}}^{\varrho}\cup (e,\varrho]\cup I_{\varrho}^{\mathfrak{e}}, \kappa \in I_{\mathfrak{e},\varrho}$ and $z\in (\varrho,1]$, then $\mathbb{U}_{[0,\varrho]}^{1}(\iota ,\mathbb{U}_{[0,\varrho]}^{1}( \kappa ,\omega ))=\mathbb{U}_{[0,\varrho]}^{1}(\iota ,1)=1=\mathbb{U}_{[0,\varrho]}^{1}(1,\omega )=\mathbb{U}_{[0,\varrho]}^{1}(\mathbb{U}_{[0,\varrho]}^{1}(\iota , \kappa ),\omega )$ and $\mathbb{U}_{[0,\varrho]}^{1}(\kappa , \mathbb{U}_{[0,\varrho]}^{1}(\iota ,\omega ))=\mathbb{U}_{[0,\varrho]}^{1}(\kappa ,1)=1$. Thus $\mathbb{U}_{[0,\varrho]}^{1}(\iota ,\mathbb{U}_{[0,\varrho]}^{1}(\kappa ,\omega))=\mathbb{U}_{[0,\varrho]}^{1}(\mathbb{U}_{[0,\varrho]}^{1}(\iota , \kappa ),\omega)$ $=\mathbb{U}_{[0,\varrho]}^{1}(\kappa ,\mathbb{U}_{[0,\varrho]}^{1}(\iota ,\omega))$.

%Combining the above cases, we obtain that $\mathbb{U}_{[0,\varrho]}^{1}(x,\mathbb{U}_{[0,\varrho]}^{1}(y,z))=\mathbb{U}_{[0,\varrho]}^{1}(\mathbb{U}_{[0,\varrho]}^{1}(x,y),z)$ for all $x,y,z\in L$ by Proposition  \ref{pro2.1}.
Therefore, $\mathbb{U}_{[0,\varrho]}^{1}$ is a uninorm on $\mathbb{L_{B}}$ with the neutral element $e$.

(2) Next we just prove  that if $I_{\varrho}^{\mathfrak{e}}\cup I_{\mathfrak{e},\varrho}\cup (\varrho,1)\neq\emptyset$, then the condition $\mathbb{U}^{*}\in \mathcal{U}_{b}$ is   necessary for that $\mathbb{U}_{[0,\varrho]}^{1}$ is a uninorm on $\mathbb{L_{B}}$.

Suppose that  $I_{\varrho}^{\mathfrak{e}}\cup I_{\mathfrak{e},\varrho}\cup (\varrho,1)\neq\emptyset$ and  $\mathbb{U}_{[0,\varrho]}^{1}(\iota , \kappa )$ is a uninorm on $\mathbb{L_{B}}$. We prove that if $\mathbb{U}^{*}(\iota , \kappa )\in[0,\mathfrak{e}]  $, then $(\iota , \kappa )\in [0,\mathfrak{e}]^{2}$. The proof can be split into all possible cases:

(i) $\mathbb{U}^{*}(\iota , \kappa )\notin[0,\mathfrak{e}]  $ for all $(\iota , \kappa )\in [0,\mathfrak{e}]\times (I_{\mathfrak{e}}^{\varrho}\cup (\mathfrak{e},\varrho])\cup (I_{\mathfrak{e}}^{\varrho}\cup (\mathfrak{e},\varrho])\times [0,\mathfrak{e}]$.

Now we just prove that  $\mathbb{U}^{*}(\iota , \kappa )\notin[0,\mathfrak{e}]  $ for all $(\iota , \kappa )\in [0,\mathfrak{e}]\times (I_{\mathfrak{e}}^{\varrho}\cup (\mathfrak{e},\varrho])$, and the other case  can be proved immediately by the commutativity property of $\mathbb{U}^{*}$. Assume that there exists $ (\iota , \kappa )\in [0,\mathfrak{e}]\times (I_{\mathfrak{e}}^{\varrho}\cup (\mathfrak{e},\varrho]) $ such that $\mathbb{U}^{*}(\iota , \kappa )\in[0,\mathfrak{e}]$. Take $z\in I_{\varrho}^{\mathfrak{e}}\cup I_{\mathfrak{e},\varrho}\cup (\varrho,1)$. Then $\mathbb{U}_{[0,\varrho]}^{1}(\iota ,\mathbb{U}_{[0,\varrho]}^{1}(\kappa ,\omega))=\mathbb{U}_{[0,\varrho]}^{1}(\iota ,1)=1$ and $\mathbb{U}_{[0,\varrho]}^{1}(\mathbb{U}_{[0,\varrho]}^{1}(\iota , \kappa ),\omega)=\mathbb{U}_{[0,\varrho]}^{1}(\mathbb{U}^{*}(\iota , \kappa ),\omega)=\omega$. Since $\omega \neq 1$, the associativity property of $\mathbb{U}_{[0,\varrho]}^{1}(\iota , \kappa )$ is violated. Thus $\mathbb{U}^{*}(\iota , \kappa )\notin[0,\mathfrak{e}]  $ for all $(\iota , \kappa )\in [0,\mathfrak{e}]\times (I_{\mathfrak{e}}^{\varrho}\cup (\mathfrak{e},\varrho])\cup (I_{\mathfrak{e}}^{\varrho}\cup (\mathfrak{e},\varrho])\times [0,\mathfrak{e}]$.

(ii) $\mathbb{U}^{*}(\iota , \kappa )\notin[0,\mathfrak{e}]  $ for all $(\iota , \kappa )\in I_{\mathfrak{e}}^{\varrho}\times I_{\mathfrak{e}}^{\varrho}$.

Assume that there exists $ (\iota , \kappa )\in I_{\mathfrak{e}}^{\varrho}\times I_{\mathfrak{e}}^{\varrho}$ such that $\mathbb{U}^{*}(\iota , \kappa )\in[0,\mathfrak{e}]$. Take $\omega \in I_{\varrho}^{\mathfrak{e}}\cup I_{\mathfrak{e},\varrho}\cup (\varrho,1)$. Then $\mathbb{U}_{[0,\varrho]}^{1}(\iota ,\mathbb{U}_{[0,\varrho]}^{1}(\kappa ,\omega))=\mathbb{U}_{[0,\varrho]}^{1}(\iota ,1)=1$ and $\mathbb{U}_{[0,\varrho]}^{1}(\mathbb{U}_{[0,\varrho]}^{1}(\iota , \kappa ),\omega)=\mathbb{U}_{[0,\varrho]}^{1}(\mathbb{U}^{*}(\iota , \kappa ),\omega)=\omega$. Since $\omega \neq 1$, the associativity property of $\mathbb{U}_{[0,\varrho]}^{1}(\iota , \kappa )$ is violated. Thus $\mathbb{U}^{*}(\iota , \kappa )\notin[0,\mathfrak{e}]  $ for all $(\iota , \kappa )\in I_{\mathfrak{e}}^{\varrho}\times I_{\mathfrak{e}}^{\varrho}$.

(iii) $\mathbb{U}^{*}(\iota , \kappa )\notin[0,\mathfrak{e}]  $ for all $(\iota , \kappa )\in (\mathfrak{e},\varrho]^{2}\cup (\mathfrak{e},\varrho]\times I_{\mathfrak{e}}^{\varrho}\cup I_{\mathfrak{e}}^{\varrho}\times (\mathfrak{e},\varrho]$.

Now we just prove that  $\mathbb{U}^{*}(\iota , \kappa )\notin[0,\mathfrak{e}]  $ for all $(\iota ,\kappa )\in (\mathfrak{e},\varrho]^{2}\cup (\mathfrak{e},\varrho]\times I_{\mathfrak{e}}^{\varrho}$, and the other case  can be proved immediately by the commutativity property of $\mathbb{U}^{*}$. By the increasingness property of $\mathbb{U}^{*}$,  we can obtain that $\kappa =\mathbb{U}^{*}(\mathfrak{e}, \kappa )\leq \mathbb{U}^{*}(\iota , \kappa )$. Since $\kappa \in I_{\mathfrak{e}}^{\varrho} \cup (\mathfrak{e},\varrho]$, we can obtain that $\mathbb{U}^{*}(\iota , \kappa )\notin [0,\mathfrak{e}]$. Thus $\mathbb{U}^{*}(\iota , \kappa )\notin[0,\mathfrak{e}]  $ for all $(\iota , \kappa )\in (\mathfrak{e},\varrho]^{2}\cup (\mathfrak{e},\varrho]\times I_{\mathfrak{e}}^{\varrho}\cup I_{\mathfrak{e}}^{\varrho}\times (\mathfrak{e},\varrho]$.

Hence, $\mathbb{U}^{*}(\iota , \kappa )\in[0,\mathfrak{e}]  $ implies $(\iota , \kappa )\in [0,\mathfrak{e}]^{2}$.
\end{proof}

\begin{remark}
In Theorem \ref{th31}, whether $q\in (0,\mathfrak{e})$ or $q\in I_{\mathfrak{e}}^{\varrho}$, if $q\nparallel \iota $ for $\iota \in I_{\mathfrak{e},\varrho}$, i.e. $q<\iota $, then $\iota \vee \kappa \vee q=\iota \vee \kappa $ for $\kappa \in I_{\mathfrak{e},\varrho}$. Moreover, if $\iota \parallel q$ for $\iota \in I_{\mathfrak{e},\varrho}$, then $\iota \vee \kappa \vee q=1$ for $\kappa \in I_{\mathfrak{e},\varrho}$. Therefore, whether $q\in (0,\mathfrak{e})$ or $q\in I_{\mathfrak{e}}^{\varrho}$ do not affect the conclusion.
\end{remark}

If we take $\mathfrak{e}=0$ in Theorem \ref{th31}, then we can obtain the existing result.

\begin{remark}\label{re311}
In Theorem \ref{th31}, if taking $\mathfrak{e}=0$, then $\mathbb{U}^{*}$ is a $t$-conorm  and $I_{\mathfrak{e}}=I_{0}=\varnothing$. In this case, the conditions in Theorem \ref{th31} naturally hold.

By the above fact, if taking $\mathfrak{e}=0$ in Theorem \ref{th31}, then the $t$-conorm $\mathbb{U}_{[0,\varrho]}^{1}:\mathbb{L_{B}}^{2}\rightarrow \mathbb{L_{B}}$ can be obtained as follows:

$\mathbb{U}_{[0,\varrho]}^{1}(\iota , \kappa )=\begin{cases}
 \mathbb{U}^{*}(\iota , \kappa ) &\mbox{if } (\iota ,\kappa )\in [0,\varrho]^{2},\\
 \iota \vee \kappa  &\mbox{if } 0\in \{\iota ,\kappa \},\\
 1 &\mbox{} otherwise.
\end{cases}$\

Obviously, $\mathbb{U}_{[0,\varrho]}^{1}$ is just the $t$-conorm $\mathbb{S}$ in Theorem \ref{th021}.
% the same as
\end{remark}

The next example illustrates the  method of uninorms  with $q\in (0,\mathfrak{e})$ in Theorem \ref{th31}.

\begin{example}
Given a  lattice $\mathbb{L_{B}}_{11}=\{0,q,\mathfrak{e},k,c,\varrho,m,t,s,d,1\}$ depicted in Fig.1.1 and a uninorm $\mathbb{U}^{*}:[0,\varrho]^{2}\rightarrow[0,\varrho]$ shown in Table \ref{Tab:01}. It is easy to see that $\mathbb{L_{B}}_{11}$ and $\mathbb{U}^{*}$ satisfy the conditions in Theorem \ref{th31} with $q\in (0,\mathfrak{e})$. Based on Theorem \ref{th31}, the uninorm $\mathbb{U}_{[0,\varrho]}^{1}:\mathbb{L_{B}}_{11}^{2}\rightarrow \mathbb{L_{B}}_{11}$ with the neutral element $\mathfrak{e}$ is defined as in Table \ref{Tab:02}.
\end{example}

\begin{minipage}{11pc}
\setlength{\unitlength}{0.75pt}\begin{picture}(600,240)
\put(266,29){$\bullet$}\put(267,22){\makebox(0,0)[l]{\footnotesize$0$}}
\put(266,60){$\bullet$}\put(257,65){\makebox(0,0)[l]{\footnotesize$q$}}
\put(266,92){$\bullet$}\put(257,96){\makebox(0,0)[l]{\footnotesize$\mathfrak{e}$}}
\put(266,124){$\bullet$}\put(257,128){\makebox(0,0)[l]{\footnotesize$c$}}
\put(266,156){$\bullet$}\put(257,160){\makebox(0,0)[l]{\footnotesize$\varrho$}}
\put(266,188){$\bullet$}\put(277,195){\makebox(0,0)[l]{\footnotesize$d$}}
\put(266,220){$\bullet$}\put(267,235){\makebox(0,0)[l]{\footnotesize$1$}}

\put(221,139){$\bullet$}\put(218,150){\makebox(0,0)[l]{\footnotesize$s$}}

\put(185,142){$\bullet$}\put(174,140){\makebox(0,0)[l]{\footnotesize$t$}}

\put(297,92){$\bullet$}\put(308,96){\makebox(0,0)[l]{\footnotesize$k$}}

\put(344,141){$\bullet$}\put(353,144){\makebox(0,0)[l]{\footnotesize$m$}}

\put(270,32){\line(0,1){30}}
\put(270,65){\line(0,1){30}}
\put(270,97){\line(0,1){30}}
\put(270,129){\line(0,1){30}}
\put(270,161){\line(0,1){30}}
\put(270,193){\line(0,1){30}}
\put(270,65){\line(-1,1){83}}
\put(268,224){\line(-1,-1){80}}
\put(270,97){\line(-1,1){45}}
\put(270,190){\line(-1,-1){45}}
\put(270,65){\line(1,1){79}}
\put(300,97){\line(-1,1){30}}
\put(270,224){\line(1,-1){79}}

\put(200,0){\emph{Fig.1.1. The lattice $\mathbb{L_{B}}_{11}$}}
\end{picture}
\end{minipage}\\

\begin{table}[htbp]
\centering
\caption{The uninorm $\mathbb{U}^{*}$ on $[0,\varrho]$.}
\label{Tab:01}

\begin{tabular}{c|c c c c c c}
\hline
  $\mathbb{U}^{*}$ & $0$ & $q$ & $\mathfrak{e}$ & $k$ & $c$ & $\varrho$ \\
\hline
  $0$ & $0$ & $0$ & $0$ & $k$ & $c$ & $\varrho$ \\

  $q$ & $0$ & $q$ & $q$ & $k$ & $c$ & $\varrho$ \\

  $\mathfrak{e}$ & $0$ & $q$ & $\mathfrak{e}$ & $k$ & $c$ & $\varrho$ \\

  $k$ & $k$ & $k$ & $k$ & $k$ & $c$ & $\varrho$ \\

  $c$ & $c$ & $c$ & $c$ & $c$ & $c$ & $\varrho$ \\

  $\varrho$ & $\varrho$ & $\varrho$ & $\varrho$ & $\varrho$ & $\varrho$ & $\varrho$ \\
\hline
\end{tabular}
\end{table}
\vspace{-0.5cm}

\begin{table}[htbp]
\centering
\caption{The uninorm $\mathbb{U}_{[0,\varrho]}^{1} $ on $\mathbb{L_{B}}_{11}$.}
\label{Tab:02}

\begin{tabular}{c|c c c c c c c c c c c c c}
\hline
  $U_{11}$ & $0$ & $q$ & $\mathfrak{e}$ & $k$ & $c$ & $\varrho$ & $m$ & $t$  & $s$ & $d$ & $1$ \\
\hline
  $0$ & $0$ & $0$ & $0$ & $k$ & $c$ & $\varrho$ & $m$ & $t$ & $s$ & $d$ & $1$ \\

  $q$ & $0$ & $q$ & $q$ & $k$ & $c$ & $\varrho$ & $m$ & $t$ & $s$ & $d$ & $1$ \\

  $\mathfrak{e}$ & $0$ & $q$ & $\mathfrak{e}$ & $k$ & $c$ & $\varrho$ & $m$ & $t$ & $s$ & $d$ & $1$ \\

  $k$ & $k$ & $k$ & $k$ & $k$ & $c$ & $\varrho$ & $1$ & $1$ & $1$ & $1$ & $1$ \\

  $c$ & $c$ & $c$ & $c$ & $c$ & $c$ & $\varrho$ & $1$ & $1$ & $1$ & $1$ & $1$ \\

  $\varrho$ & $\varrho$ & $\varrho$ & $\varrho$ & $\varrho$ & $\varrho$ & $\varrho$ & $1$ & $1$ & $1$ & $1$ & $1$ \\

  $m$ & $m$ & $m$ & $m$ & $1$ & $1$ & $1$ & $m$ & $1$ & $1$ & $1$ & $1$ \\

  $t$ & $t$ & $t$ & $t$ & $1$ & $1$ & $1$ & $1$ & $t$ & $1$ & $1$ & $1$ \\

  $s$ & $s$ & $s$ & $s$ & $1$ & $1$ & $1$ & $1$ & $1$ & $1$ & $1$ & $1$ \\

  $d$ & $d$ & $d$ & $d$ & $1$ & $1$ & $1$ & $1$ & $1$ & $1$ & $1$ & $1$ \\

  $1$ & $1$ & $1$ & $1$ & $1$ & $1$ & $1$ & $1$ & $1$ & $1$ & $1$ & $1$ \\
\hline
\end{tabular}
\end{table}\

Moreover, we give the  example to illustrate the construction method of uninorms  in Theorem \ref{th31}.

\begin{example}
Given a lattice $\mathbb{L_{B}}_{12}=\{0,f,\mathfrak{e},c,q,\varrho,m,t,s,d,1\}$ depicted in Fig.1.2 and a uninorm $\mathbb{U}^{*}:[0,\varrho]^{2}\rightarrow[0,\varrho]$ shown in Table \ref{Tab:03}. It is easy to see that $\mathbb{L_{B}}_{12}$ and $\mathbb{U}^{*}$ satisfy the conditions in Theorem \ref{th31} with $q\in I_{\mathfrak{e}}^{\varrho}$. Based on Theorem \ref{th31}, the uninorm $\mathbb{U}_{[0,\varrho]}^{1}:\mathbb{L_{B}}_{12}^{2}\rightarrow \mathbb{L_{B}}_{12}$
 %with the neutral element $\mathfrak{e}$
 is defined as in Table \ref{Tab:04}.
\end{example}

\begin{minipage}{11pc}
\setlength{\unitlength}{0.75pt}\begin{picture}(600,240)
\put(266,29){$\bullet$}\put(267,22){\makebox(0,0)[l]{\footnotesize$0$}}
\put(266,60){$\bullet$}\put(257,59){\makebox(0,0)[l]{\footnotesize$f$}}
\put(266,92){$\bullet$}\put(257,96){\makebox(0,0)[l]{\footnotesize$\mathfrak{e}$}}
\put(266,124){$\bullet$}\put(257,128){\makebox(0,0)[l]{\footnotesize$c$}}
\put(266,156){$\bullet$}\put(257,160){\makebox(0,0)[l]{\footnotesize$\varrho$}}
\put(266,188){$\bullet$}\put(277,195){\makebox(0,0)[l]{\footnotesize$d$}}
\put(266,220){$\bullet$}\put(267,235){\makebox(0,0)[l]{\footnotesize$1$}}
\put(219,106){$\bullet$}\put(218,120){\makebox(0,0)[l]{\footnotesize$q$}}
\put(171,122){$\bullet$}\put(170,140){\makebox(0,0)[l]{\footnotesize$t$}}
\put(314,139){$\bullet$}\put(305,145){\makebox(0,0)[l]{\footnotesize$k$}}
\put(344,141){$\bullet$}\put(353,144){\makebox(0,0)[l]{\footnotesize$m$}}

\put(270,32){\line(0,1){30}}
\put(270,65){\line(0,1){30}}
\put(270,97){\line(0,1){30}}
\put(270,129){\line(0,1){30}}
\put(270,161){\line(0,1){30}}
\put(270,193){\line(0,1){30}}
\put(270,31){\line(-1,1){95}}
\put(270,224){\line(-1,-1){95}}

\put(270,63){\line(-1,1){47}}
\put(270,158){\line(-1,-1){47}}

\put(270,95){\line(1,1){47}}
\put(270,192){\line(1,-1){47}}

\put(270,65){\line(1,1){79}}
\put(270,224){\line(1,-1){79}}

\put(200,0){\emph{Fig.1.2. The lattice $\mathbb{L_{B}}_{12}$}}
\end{picture}
\end{minipage}\

\begin{table}[htbp]
\centering
\caption{The uninorm $\mathbb{U}^{*}$ on $[0,\varrho]$.}
\label{Tab:03}\

\begin{tabular}{c|c c c c c c}
\hline
  $\mathbb{U}^{*}$ & $0$ & $f$ & $\mathfrak{e}$ & $c$ & $q$ & $\varrho$ \\
\hline
  $0$ & $0$ & $0$ & $0$ & $c$ & $q$ & $\varrho$ \\

  $f$ & $0$ & $f$ & $f$ & $c$ & $q$ & $\varrho$ \\

  $\mathfrak{e}$ & $0$ & $f$ & $\mathfrak{e}$ & $c$ & $q$ & $\varrho$ \\

  $c$ & $c$ & $c$ & $c$ & $c$ & $q$ & $\varrho$ \\

  $q$ & $q$ & $q$ & $q$ & $q$ & $q$ & $\varrho$ \\

  $\varrho$ & $\varrho$ & $\varrho$ & $\varrho$ & $\varrho$ & $\varrho$ & $\varrho$ \\
\hline
\end{tabular}
\end{table}
\vspace{-0.5cm}

\begin{table}[htbp]
\centering
\caption{The uninorm $\mathbb{U}_{[0,\varrho]}^{1}$ on $\mathbb{L_{B}}_{12}$.}
\label{Tab:04}\

\begin{tabular}{c|c c c c c c c c c c c c c}
\hline
  $U_{12}$ & $0$ & $f$ & $\mathfrak{e}$ & $c$ & $q$ & $\varrho$ & $s$ & $t$  & $m$ & $d$ & $1$ \\
\hline
  $0$ & $0$ & $0$ & $0$ & $c$ & $q$ & $\varrho$ & $s$ & $t$ & $m$ & $d$ & $1$ \\

  $f$ & $0$ & $f$ & $f$ & $c$ & $q$ & $\varrho$ & $s$ & $t$ & $m$ & $d$ & $1$ \\

  $\mathfrak{e}$ & $0$ & $q$ & $\mathfrak{e}$ & $c$ & $q$ & $\varrho$ & $s$ & $t$ & $m$ & $d$ & $1$ \\

  $c$ & $c$ & $c$ & $c$ & $c$ & $q$ & $\varrho$ & $1$ & $1$ & $1$ & $1$ & $1$ \\

  $q$ & $q$ & $q$ & $q$ & $q$ & $q$ & $\varrho$ & $1$ & $1$ & $1$ & $1$ & $1$ \\

  $\varrho$ & $\varrho$ & $\varrho$ & $\varrho$ & $\varrho$ & $\varrho$ & $\varrho$ & $1$ & $1$ & $1$ & $1$ & $1$ \\

  $s$ & $s$ & $s$ & $s$ & $1$ & $1$ & $1$ & $1$ & $1$ & $1$ & $1$ & $1$ \\

  $t$ & $t$ & $t$ & $t$ & $1$ & $1$ & $1$ & $1$ & $1$ & $1$ & $1$ & $1$ \\

  $m$ & $m$ & $m$ & $m$ & $1$ & $1$ & $1$ & $1$ & $1$ & $1$ & $1$ & $1$ \\

  $d$ & $d$ & $d$ & $d$ & $1$ & $1$ & $1$ & $1$ & $1$ & $1$ & $1$ & $1$ \\

  $1$ & $1$ & $1$ & $1$ & $1$ & $1$ & $1$ & $1$ & $1$ & $1$ & $1$ & $1$ \\
\hline
\end{tabular}\
\
\
\end{table}

\begin{remark}\label{re312}
In Theorem \ref{th31},  the conditions $\iota \vee \kappa =1$ for all $\iota ,\kappa \in I_{\mathfrak{e},\varrho}$ with $\iota \neq \kappa $, and $\iota \vee q=1$ for all $\iota \in I_{\mathfrak{e}}^{\varrho} \cap I_{q}$ can not be omitted, in general.
\end{remark}

The next example illustrates the fact in Remark \ref{re312}.
%That is, if the conditions in Theorem \ref{th31} do not hold, then the increasingness of $\mathbb{U}_{[0,\varrho]}^{1}$ can be violated.
%\newpage

\begin{example}\label{ex312}
Given a lattice $\mathbb{L_{B}}_{13}=\{0,q,\mathfrak{e},\varrho,m,t,s,d,1\}$ depicted in Fig.1.3 and a uninorm $\mathbb{U}^{*}:[0,\varrho]^{2}\rightarrow[0,\varrho]$ shown in Table \ref{Tab:0121}. We can see that $\mathbb{U}^{*}\in \mathcal{U}_{b}$ and $\iota \parallel \kappa $ for all $\iota \in I_{\mathfrak{e},\varrho}\cap I^{q}$ and $y\in I_{\mathfrak{e}}^{\varrho} $. Since $t\vee m =d<1$ and $m\vee q=d<1$, the conditions that  $\iota \vee \kappa =1$ for all $\iota ,\kappa \in I_{\mathfrak{e},\varrho}$ with $\iota \neq \kappa $ and $\iota \vee q=1$ for all $\iota \in I_{\mathfrak{e},\varrho} \cap I_{q}$ in Theorem \ref{th31} do not hold. Based on  Theorem \ref{th31} with $q\in (0,\mathfrak{e})$, we can obtain a function $\mathbb{U}_{[0,\varrho]}^{1}$ on $\mathbb{L_{B}}_{13}$, shown in Table \ref{Tab:0122}. Since $\mathbb{U}_{[0,\varrho]}^{1}(t,m)=d<1$ and $\mathbb{U}_{[0,\varrho]}^{1}(s,m)=1$,  $\mathbb{U}_{[0,\varrho]}^{1}$ does not satisfy the increasingness. Thus, $\mathbb{U}_{[0,\varrho]}^{1}$ is not a uninorm on $\mathbb{L_{B}}_{13}$.
\end{example}

\begin{minipage}{11pc}
\setlength{\unitlength}{0.75pt}\begin{picture}(600,200)
\put(266,29){$\bullet$}\put(267,22){\makebox(0,0)[l]{\footnotesize$0$}}
\put(266,60){$\bullet$}\put(257,65){\makebox(0,0)[l]{\footnotesize$q$}}
\put(266,92){$\bullet$}\put(257,96){\makebox(0,0)[l]{\footnotesize$\mathfrak{e}$}}

\put(266,126){$\bullet$}\put(257,130){\makebox(0,0)[l]{\footnotesize$\varrho$}}
\put(266,158){$\bullet$}\put(277,165){\makebox(0,0)[l]{\footnotesize$d$}}
\put(266,190){$\bullet$}\put(267,205){\makebox(0,0)[l]{\footnotesize$1$}}

\put(217,77){$\bullet$}\put(217,68){\makebox(0,0)[l]{\footnotesize$s$}}

\put(202,93){$\bullet$}\put(191,95){\makebox(0,0)[l]{\footnotesize$t$}}

\put(330,93){$\bullet$}\put(341,95){\makebox(0,0)[l]{\footnotesize$m$}}

\put(270,32){\line(0,1){170}}

\put(270,32){\line(-1,1){65}}
\put(270,163){\line(-1,-1){65}}

\put(270,131){\line(-1,-1){50}}

\put(270,32){\line(1,1){65}}

\put(270,164){\line(1,-1){65}}

\put(200,0){\emph{Fig.1.3. The lattice $\mathbb{L_{B}}_{13}$}}
\end{picture}
\end{minipage}\

\begin{table}[htbp]
\centering
\caption{The uninorm $\mathbb{U}^{*}$ on $[0,\varrho]$.}
\label{Tab:0121}\

\begin{tabular}{c|c c c c c c c}
\hline
  $\mathbb{U}^{*}$ & $0$ & $q$ & $\mathfrak{e}$ & $\varrho$ \\
\hline
  $0$ & $0$ & $0$ & $0$ & $\varrho$ \\

  $q$ & $0$ & $q$ & $q$ & $\varrho$ \\

  $\mathfrak{e}$ & $0$ & $q$ & $\mathfrak{e}$ & $\varrho$ \\

  $\varrho$ & $\varrho$ & $\varrho$ & $\varrho$ & $\varrho$ \\
\hline
\end{tabular}
\end{table}

\begin{table}[htbp]
\centering
\caption{The function $\mathbb{U}_{[0,\varrho]}^{1} $ on $\mathbb{L_{B}}_{13}$.}
\label{Tab:0122}\

\begin{tabular}{c|c c c c c c c c c}
\hline
  $U_{13}$ & $0$ & $q$ & $\mathfrak{e}$ & $\varrho$ & $s$ & $t$ & $m$ & $d$ & $1$ \\
\hline
  $0$ & $0$ & $0$ & $0$ & $\varrho$ & $s$ & $t$ & $m$ & $d$ & $1$ \\

  $q$ & $0$ & $q$ & $q$ & $\varrho$ & $s$ & $t$ & $m$ & $1$ & $1$ \\

  $\mathfrak{e}$ & $0$ & $q$ & $\mathfrak{e}$ & $\varrho$ & $s$ & $t$ & $m$ & $d$ & $1$ \\

  $\varrho$ & $\varrho$ & $\varrho$ & $\varrho$ & $\varrho$ & $1$ & $1$ & $1$ & $1$ & $1$ \\

  $s$ & $s$ & $s$ & $s$ & $1$ & $1$ & $1$ & $1$ & $1$ & $1$ \\

  $t$ & $t$ & $t$ & $t$ & $1$ & $1$ & $d$ & $d$ & $1$ & $1$ \\

  $m$ & $m$ & $m$ & $m$ & $1$ & $1$ & $d$ & $d$ & $1$ & $1$ \\

  $d$ & $d$ & $d$ & $d$ & $1$ & $1$ & $1$ & $1$ & $1$ & $1$ \\

  $1$ & $1$ & $1$ & $1$ & $1$ & $1$ & $1$ & $1$ & $1$ & $1$ \\
\hline
\end{tabular}
\end{table}

\begin{remark}\label{re32}
% Let $\mathbb{U}_{[0,\varrho]}^{1}$ be a uninorm in  Theorem \ref{th31}.\\
$(1)$ If $\varrho=1$, then $\mathbb{U}_{[0,\varrho]}^{1}=\mathbb{U}^{*}$.\\
$(2)$ If $\mathbb{U}^{*}\in \mathcal{U}_{b}$, then $\mathbb{U}_{[0,\varrho]}^{1}\in \mathcal{U}_{b}$.\\
$(3)$ $\mathbb{U}_{[0,\varrho]}^{1}\in \mathcal{U}_{max}$ if and only if $\mathbb{U}^{*}\in \mathcal{U}_{max}$.
\end{remark}

Next, we illustrate that the function $\mathbb{U}$ given by $(1)$ with $q \in I_{\varrho}^{\mathfrak{e}}$ can be a uninorm
on bounded lattices under some conditions.

\begin{theorem}\label{th33}
Let $\varrho \in \mathbb{L_{B}} \setminus\{0,1\}$, $q\in I_{\varrho}^{\mathfrak{e}}$, $\mathbb{U}^{*}$ be a uninorm on $[0,\varrho]$ with a neutral element $e$ and $\mathbb{U}_{[0,\varrho]}^{2}$ be a function given by $(\ref{eq1})$. Suppose that $\iota \vee q =1$ for all $\iota \in I_{\mathfrak{e},\varrho}\cap I_{q}$.

$(1)$ Let us assume that $\mathbb{U}^{*}\in \mathcal{U}_{b}$. Then $\mathbb{U}_{[0,\varrho]}^{2}$ is a uninorm on $\mathbb{L_{B}}$ with the neutral element $\mathfrak{e} \in \mathbb{L_{B}} $ if and only if $\iota \parallel y$ for all $\iota \in I_{\mathfrak{e},\varrho}\cap I^{q}$ and $\kappa \in I_{\mathfrak{e}}^{\varrho} $.

$(2)$ Moreover, let us assume that $I_{\varrho}^{\mathfrak{e}}\cup I_{\mathfrak{e},\varrho}\cup (\varrho,1)\neq\emptyset$. Then $\mathbb{U}_{[0,\varrho]}^{2}$ is a uninorm on $\mathbb{L_{B}}$ with the neutral element $\mathfrak{e} \in \mathbb{L_{B}} $ if and only if $\mathbb{U}^{*}\in \mathcal{U}_{b}$ and $\iota \parallel \kappa$ for all $\iota \in I_{\mathfrak{e},\varrho}\cap I^{q}$ and $\kappa \in I_{\mathfrak{e}}^{\varrho} $.

\end{theorem}

 \begin{proof}
 (1) Necessity. Let $\mathbb{U}_{[0,\varrho]}^{2}$ be a uninorm on $\mathbb{L_{B}}$ with a neutral element $\mathfrak{e}$. We need to prove that $\iota \parallel \kappa$ for all  $\iota \in I_{\mathfrak{e},\varrho}\cap I^{q}$ and $\kappa \in I_{\mathfrak{e}}^{\varrho} $.

Assume that there exist $\iota \in I_{\mathfrak{e},\varrho} \cap I^{q}$ and $y\in I_{\mathfrak{e}}^{\varrho} $  such that $\iota \nparallel \kappa $, i.e., $\kappa < \iota $ and $\iota <q$. Then $\mathbb{U}_{[0,\varrho]}^{2}(\iota , \kappa )=1$ and $\mathbb{U}_{[0,\varrho]}^{2}(\iota ,\iota )=\iota \vee q=q$. Since $q<1$, the  increasingness property of $\mathbb{U}_{[0,\varrho]}^{2}$ is violated. Thus $\iota \parallel \kappa $ for all $\iota \in I_{\mathfrak{e},\varrho}\cap I^{q} $ and $\kappa \in I_{\mathfrak{e}}^{\varrho} $.

Sufficiency. It is obvious that  $\mathbb{U}_{[0,\varrho]}^{2}$ is commutative  and  $\mathfrak{e}$  is the neutral element of $\mathbb{U}_{[0,\varrho]}^{2}$. Thus, we just prove the increasingness and the associativity of $\mathbb{U}_{[0,\varrho]}^{2}$.

I. Increasingness: It can be obtained by the proof  of  Theorem \ref{th31} in a similar way.

II. Associativity:  By Proposition \ref{pro2.1} and  Theorem \ref{th31}, it is enough to check only those cases that are different from the cases in Theorem \ref{th31}.

1. Suppose that $\iota , \kappa , \omega \in I_{\mathfrak{e},\varrho}$.

1.1. If $\iota , \kappa ,\omega \nparallel q$, then $\mathbb{U}_{[0,\varrho]}^{2}(\iota ,\mathbb{U}_{[0,\varrho]}^{2}(\kappa ,\omega))=\mathbb{U}_{[0,\varrho]}^{2}(\iota ,\kappa \vee \omega \vee q)=\mathbb{U}_{[0,\varrho]}^{2}(\iota ,q)=1=\mathbb{U}_{[0,\varrho]}^{2}(q,\omega)=\mathbb{U}_{[0,\varrho]}^{2}(\iota \vee \kappa \vee q,\omega)= \mathbb{U}_{[0,\varrho]}^{2}(\mathbb{U}_{[0,\varrho]}^{2}(\iota , \kappa ),\omega)$  and $\mathbb{U}_{[0,\varrho]}^{2}(\kappa ,\mathbb{U}_{[0,\varrho]}^{2}(\iota ,\omega))=\mathbb{U}_{[0,\varrho]}^{2}(\kappa , \iota \vee \omega \vee q)=\mathbb{U}_{[0,\varrho]}^{2}(\kappa ,q)=1$.

1.2 Assume that there exists $ \iota \in I_{\mathfrak{e},\varrho}$ such that $\iota \parallel q$.

1.2.1. If $\iota \nparallel q$ and $\kappa ,\omega \parallel q$, then $\mathbb{U}_{[0,\varrho]}^{2}(\iota ,\mathbb{U}_{[0,\varrho]}^{2}(\kappa ,\omega ))=\mathbb{U}_{[0,\varrho]}^{2}(\iota ,\kappa \vee \omega \vee q)=\mathbb{U}_{[0,\varrho]}^{2}(\iota ,1)=1=\mathbb{U}_{[0,\varrho]}^{2}(1,\omega)=\mathbb{U}_{[0,\varrho]}^{2}(\iota \vee \kappa \vee q,\omega)= \mathbb{U}_{[0,\varrho]}^{2}(\mathbb{U}_{[0,\varrho]}^{2}(\iota ,\kappa ),\omega)$ and $\mathbb{U}_{[0,\varrho]}^{2}(\kappa ,\mathbb{U}_{[0,\varrho]}^{2}(\iota ,\omega))=\mathbb{U}_{[0,\varrho]}^{2}(\kappa ,\iota \vee \omega \vee q)=\mathbb{U}_{[0,\varrho]}^{2}(\kappa ,1)=1$.

1.2.2. If $\iota ,\kappa \nparallel q$ and $\omega \parallel q$, then $\mathbb{U}_{[0,\varrho]}^{2}(\iota ,\mathbb{U}_{[0,\varrho]}^{2}(\kappa , \omega ))=\mathbb{U}_{[0,\varrho]}^{2}(\iota ,1)=1=\mathbb{U}_{[0,\varrho]}^{2}(\kappa ,\iota \vee \omega \vee q)= \mathbb{U}_{[0,\varrho]}^{2}(\kappa ,\mathbb{U}_{[0,\varrho]}^{2}(\iota ,\omega ))$ and $\mathbb{U}_{[0,\varrho]}^{2}(\mathbb{U}_{[0,\varrho]}^{2}(\iota , \kappa ), \omega)=\mathbb{U}_{[0,\varrho]}^{2}(\iota \vee \kappa \vee q, \omega)=\mathbb{U}_{[0,\varrho]}^{2}(q, \omega)=1$.

1.2.3. If $\iota , \kappa , \omega \parallel q$, then $\mathbb{U}_{[0,\varrho]}^{2}(\iota ,\mathbb{U}_{[0,\varrho]}^{2}(\kappa , \omega))=\mathbb{U}_{[0,\varrho]}^{2}(\iota ,1)=1=\mathbb{U}_{[0,\varrho]}^{2}(1, \omega)=\mathbb{U}_{[0,\varrho]}^{2}(\iota \vee \kappa \vee q, \omega)= \mathbb{U}_{[0,\varrho]}^{2}(\mathbb{U}_{[0,\varrho]}^{2}(\iota , \kappa ), \omega)$ and $\mathbb{U}_{[0,\varrho]}^{2}(\kappa ,\mathbb{U}_{[0,\varrho]}^{2}(\iota , \omega))=\mathbb{U}_{[0,\varrho]}^{2}(\kappa ,1)=1$.

(2) It can be proved with the proof of Theorem \ref{th31}(2) in a similar way.
\end{proof}\

If we take $e=0$ in Theorem \ref{th33}, then we can obtain the existing result in the literature.

\begin{remark}\label{re333}
In Theorem \ref{th33}, if taking $\mathfrak{e}=0$, then $\mathbb{U}^{*}$ ia a $t$-conorm and $I_{\mathfrak{e}}=I_{0}=\varnothing$. In this case, the condition in Theorem \ref{th33} naturally holds.

By the above fact, if taking $\mathfrak{e}=0$ in Theorem \ref{th33}, then we obtain the $t$-conorm $\mathbb{U}_{[0,\varrho]}^{2}:\mathbb{L_{B}}^{2}\rightarrow \mathbb{L_{B}}$ as follows:

$\mathbb{U}_{[0,\varrho]}^{2}(\iota ,\kappa )=\begin{cases}
 \mathbb{U}^{*}(\iota , \kappa ) &\mbox{if } (\iota , \kappa )\in [0,\varrho]^{2},\\
 \iota \vee \kappa &\mbox{if } 0\in \{\iota , \kappa \},\\
 1 &\mbox{} otherwise.
\end{cases}$\

Obviously, $\mathbb{U}_{[0,\varrho]}^{2}$ is just the $t$-conorm $\mathbb{S}$ in Theorem \ref{th021}.

\end{remark}

The next example illustrates the  method of uninorms  in Theorem \ref{th33}.

\begin{example}
Given a  lattice $\mathbb{L_{B}}_{21}=\{0,f,\mathfrak{e},c,\varrho,q,m,t,d,1\}$ depicted in Fig.2.1 and a uninorm $\mathbb{U}^{*}:[0,\varrho]^{2}\rightarrow[0,\varrho]$ shown in Table \ref{Tab:05}. We can see that $\mathbb{L_{B}}_{21}$ and $\mathbb{U}^{*}$ satisfy the conditions in Theorem \ref{th33}. By  Theorem \ref{th33}, the uninorm  $\mathbb{U}_{[0,\varrho]}^{2}:\mathbb{L_{B}}_{21}^{2}\rightarrow \mathbb{L_{B}}_{21}$  is defined as in Table \ref{Tab:06}.
\end{example}

\begin{minipage}{11pc}
\setlength{\unitlength}{0.75pt}\begin{picture}(600,230)
\put(266,29){$\bullet$}\put(267,22){\makebox(0,0)[l]{\footnotesize$0$}}
\put(266,60){$\bullet$}\put(257,59){\makebox(0,0)[l]{\footnotesize$f$}}
\put(266,92){$\bullet$}\put(257,96){\makebox(0,0)[l]{\footnotesize$\mathfrak{e}$}}
\put(266,124){$\bullet$}\put(257,128){\makebox(0,0)[l]{\footnotesize$c$}}
\put(266,156){$\bullet$}\put(257,160){\makebox(0,0)[l]{\footnotesize$\varrho$}}
\put(266,188){$\bullet$}\put(277,195){\makebox(0,0)[l]{\footnotesize$d$}}
\put(266,220){$\bullet$}\put(267,235){\makebox(0,0)[l]{\footnotesize$1$}}
\put(219,106){$\bullet$}\put(218,120){\makebox(0,0)[l]{\footnotesize$m$}}
\put(171,122){$\bullet$}\put(170,140){\makebox(0,0)[l]{\footnotesize$t$}}
\put(314,139){$\bullet$}\put(305,145){\makebox(0,0)[l]{\footnotesize$q$}}

\put(270,32){\line(0,1){30}}
\put(270,65){\line(0,1){30}}
\put(270,97){\line(0,1){30}}
\put(270,129){\line(0,1){30}}
\put(270,161){\line(0,1){30}}
\put(270,193){\line(0,1){30}}

\put(270,31){\line(-1,1){95}}
\put(270,224){\line(-1,-1){95}}

\put(270,63){\line(-1,1){47}}
\put(270,158){\line(-1,-1){47}}

\put(270,95){\line(1,1){47}}
\put(270,192){\line(1,-1){47}}

\put(200,0){\emph{Fig.2.1. The lattice $\mathbb{L_{B}}_{21}$}}
\end{picture}
\end{minipage}

\begin{table}[htbp]
\centering
\caption{The uninorm $\mathbb{U}^{*}$ on $[0,\varrho]$.}
\label{Tab:05}\

\begin{tabular}{c|c c c c c c}
\hline
  $\mathbb{U}^{*}$ & $0$ & $f$ & $\mathfrak{e}$ & $c$  & $\varrho$ \\
\hline
  $0$ & $0$ & $0$ & $0$ & $c$ & $\varrho$ \\

  $f$ & $0$ & $f$ & $f$ & $c$ & $\varrho$ \\

  $\mathfrak{e}$ & $0$ & $f$ & $\mathfrak{e}$ & $c$ & $\varrho$ \\

  $c$ & $c$ & $c$ & $c$ & $c$ & $\varrho$ \\

  $\varrho$ & $\varrho$ & $\varrho$ & $\varrho$ & $\varrho$ & $\varrho$ \\
\hline
\end{tabular}
\end{table}
\vspace{-0.5cm}

\begin{table}[htbp]
\centering
\caption{The uninorm $\mathbb{U}_{[0,\varrho]}^{2} $ on $\mathbb{L_{B}}_{21}$.}
\label{Tab:06}\

\begin{tabular}{c|c c c c c c c c c c c c c}
\hline
  $U_{21}$ & $0$ & $f$ & $\mathfrak{e}$ & $c$  & $\varrho$ & $q$ & $t$ & $m$ & $d$ & $1$ \\
\hline
  $0$ & $0$ & $0$ & $0$ & $c$ & $\varrho$ & $q$ & $t$ & $m$ & $d$ & $1$ \\

  $f$ & $0$ & $f$ & $f$ & $c$ & $\varrho$ & $q$ & $t$ & $m$ & $d$ & $1$ \\

  $\mathfrak{e}$ & $0$ & $q$ & $\mathfrak{e}$ & $c$ & $\varrho$ & $q$ & $t$ & $m$ & $d$ & $1$ \\

  $c$ & $c$ & $c$ & $c$ & $c$ & $\varrho$ & $1$ & $1$ & $1$ & $1$ & $1$ \\

  $\varrho$ & $\varrho$ & $\varrho$ & $\varrho$ & $\varrho$ & $\varrho$ & $1$ & $1$ & $1$ & $1$ & $1$ \\

  $q$ & $q$ & $q$ & $q$ & $1$ & $1$ & $1$ & $1$ & $1$ & $1$ & $1$ \\

  $t$ & $t$ & $t$ & $t$ & $1$ & $1$ & $1$ & $1$ & $1$ & $1$ & $1$ \\

  $m$ & $m$ & $m$ & $m$ & $1$ & $1$ & $1$ & $1$ & $1$ & $1$ & $1$ \\

  $d$ & $d$ & $d$ & $d$ & $1$ & $1$ & $1$ & $1$ & $1$ & $1$ & $1$ \\

  $1$ & $1$ & $1$ & $1$ & $1$ & $1$ & $1$ & $1$ & $1$ & $1$ & $1$ \\
\hline
\end{tabular}
\end{table}

\begin{remark}\label{re332}
In Theorem \ref{th33},   the condition that $\iota \vee q=1$ for all $\iota \in I_{\mathfrak{e}}^{\varrho} \cap I_{q}$ can not be omitted, in general.
\end{remark}

The next example illustrates the fact in Remark \ref{re332}.
%That is, if the conditions in Theorem \ref{th33} do not hold, then the increasingness of $\mathbb{U}_{[0,\varrho]}^{2}$ can be violated.

\begin{example}\label{ex332}
Given a  lattice $\mathbb{L_{B}}_{22}=\{0,\mathfrak{e},\varrho,m,t,s,q,d,1\}$ depicted in Fig.2.2 and a uninorm $\mathbb{U}^{*}:[0,\varrho]^{2}\rightarrow[0,\varrho]$ shown in Table \ref{Tab:0321}. We can see that $\mathbb{U}^{*}\in \mathcal{U}_{b}$ and $\iota \parallel \kappa $ for all $\iota \in I_{\mathfrak{e},\varrho}\cap I^{q}$ and $\kappa \in I_{\mathfrak{e}}^{\varrho} $. Since $m\vee q=d<1$, the condition that $\iota \vee q=1$ for all $\iota \in I_{\mathfrak{e}}^{\varrho} \cap I_{q}$ in Theorem \ref{th33} do not hold. By  Theorem \ref{th33}, we can  obtain a function $\mathbb{U}_{[0,\varrho]}^{2}$ on $\mathbb{L_{B}}_{22}$, shown in Table \ref{Tab:0322}. Since $\mathbb{U}_{[0,\varrho]}^{2}(t,m)=d<1$ and $\mathbb{U}_{[0,\varrho]}^{2}(t,s)=1$,   $\mathbb{U}_{[0,\varrho]}^{2}$ does not satisfy the increasingness. Thus, $\mathbb{U}_{[0,\varrho]}^{2}$ is not a uninorm on $\mathbb{L_{B}}_{22}$.
\end{example}

\begin{minipage}{11pc}
\setlength{\unitlength}{0.75pt}\begin{picture}(600,200)
\put(266,29){$\bullet$}\put(267,22){\makebox(0,0)[l]{\footnotesize$0$}}

\put(266,80){$\bullet$}\put(257,76){\makebox(0,0)[l]{\footnotesize$\mathfrak{e}$}}

\put(266,130){$\bullet$}\put(257,126){\makebox(0,0)[l]{\footnotesize$\varrho$}}
\put(266,158){$\bullet$}\put(257,164){\makebox(0,0)[l]{\footnotesize$d$}}
\put(266,190){$\bullet$}\put(267,205){\makebox(0,0)[l]{\footnotesize$1$}}

\put(315,79){$\bullet$}\put(320,76){\makebox(0,0)[l]{\footnotesize$s$}}

\put(202,94){$\bullet$}\put(194,100){\makebox(0,0)[l]{\footnotesize$t$}}

\put(228,120){$\bullet$}\put(220,128){\makebox(0,0)[l]{\footnotesize$q$}}

\put(330,94){$\bullet$}\put(340,102){\makebox(0,0)[l]{\footnotesize$m$}}

\put(270,32){\line(0,1){50}}
\put(270,87){\line(0,1){50}}
\put(270,131){\line(0,1){30}}
\put(270,163){\line(0,1){30}}

\put(270,32){\line(-1,1){65}}
\put(270,163){\line(-1,-1){65}}

\put(270,133){\line(1,-1){50}}

\put(270,85){\line(-1,1){38}}

\put(270,32){\line(1,1){65}}

\put(270,163){\line(1,-1){65}}

\put(200,-10){\emph{Fig.2.2. The lattice $\mathbb{L_{B}}_{22}$}}
\end{picture}
\end{minipage}\\

\begin{table}[htbp]
\centering
\caption{The uninorm $\mathbb{U}^{*}$ on $[0,\varrho]$.}
\label{Tab:0321}\

\begin{tabular}{c| c c c}
\hline
  $\mathbb{U}^{*}$ & $0$ & $\mathfrak{e}$ & $\varrho$ \\
\hline
  $0$ & $0$ & $0$ & $\varrho$ \\

  $\mathfrak{e}$ & $0$ & $\mathfrak{e}$ & $\varrho$ \\

  $\varrho$ & $\varrho$ & $\varrho$ & $\varrho$ \\
\hline
\end{tabular}
\end{table}

\begin{table}[htbp]
\centering
\caption{The function $\mathbb{U}_{[0,\varrho]}^{2} $ on $\mathbb{L_{B}}_{22}$.}
\label{Tab:0322}\

\begin{tabular}{c|c c c c c c c c c}
\hline
  $U_{22}$ & $0$ & $\mathfrak{e}$ & $\varrho$ & $s$ & $t$ & $m$ & $q$ & $d$ & $1$ \\
\hline
  $0$ & $0$ & $0$ & $\varrho$ & $s$ & $t$ & $m$ & $q$ & $d$ & $1$ \\

  $\mathfrak{e}$ & $0$ & $\mathfrak{e}$ & $\varrho$ & $s$ & $t$ & $m$ & $q$ & $d$ & $1$ \\

  $\varrho$ & $\varrho$ & $\varrho$ & $\varrho$ & $1$ & $1$ & $1$ & $1$ & $1$ & $1$ \\

  $s$ & $s$ & $s$ & $1$ & $1$ & $1$ & $1$ & $1$ & $1$ & $1$ \\

  $t$ & $t$ & $t$ & $1$ & $1$ & $d$ & $d$ & $1$ & $1$ & $1$ \\

  $m$ & $m$ & $m$ & $1$ & $1$ & $d$ & $d$ & $1$ & $1$ & $1$ \\

  $q$ & $q$ & $q$ & $1$ & $1$ & $1$ & $1$ & $1$ & $1$ & $1$ \\

  $d$ & $d$ & $d$ & $1$ & $1$ & $1$ & $1$ & $1$ & $1$ & $1$ \\

  $1$ & $1$ & $1$ & $1$ & $1$ & $1$ & $1$ & $1$ & $1$ & $1$ \\
\hline
\end{tabular}
\end{table}

\begin{remark}\label{re33}
% Let $\mathbb{U}_{[0,\varrho]}^{2}$ be a uninorm in  Theorem \ref{th33}.\\
$(1)$ If $\varrho=1$, then $\mathbb{U}_{[0,\varrho]}^{2}=\mathbb{U}^{*}$.\\
$(2)$ If $\mathbb{U}^{*}\in \mathcal{U}_{b}$, then $\mathbb{U}_{[0,\varrho]}^{2}\in \mathcal{U}_{b}$.\\
$(3)$ $\mathbb{U}_{[0,\varrho]}^{2}\in \mathcal{U}_{max}$ if and only if $\mathbb{U}^{*}\in \mathcal{U}_{max}$.
\end{remark}

Similarly, let $\sigma \in \mathbb{L_{B}} \setminus\{0,1\}$, $p\in \mathbb{L_{B}}$ and $\mathbb{U}^{*}$ be a uninorm on $[\sigma,1]$ with a neutral element $e$. Then we define a function $\mathbb{U}:\mathbb{L_{B}} ^{2}\rightarrow \mathbb{L_{B}}$ by
\begin{flalign}\label{eq2}
\ \ \ \ &\ \mathbb{U}(\iota ,\kappa)=\begin{cases}
\mathbb{U}^{*}(\iota , \kappa ) &\mbox{if } (\iota ,\kappa )\in [\sigma,1]^{2},\\
\iota  &\mbox{if } (\iota ,\kappa )\in (\mathbb{L_{B}} \setminus[\sigma,1])\times [\mathfrak{e},1],\\
\kappa  &\mbox{if } (\iota , \kappa )\in [\mathfrak{e},1]\times (\mathbb{L_{B}} \setminus[\sigma,1]),\\
\iota \wedge \kappa \wedge p &\mbox{if } (\iota , \kappa )\in I_{\mathfrak{e},\sigma}\times I_{\mathfrak{e},\sigma},\\
0 &\mbox{}otherwise.\\
\end{cases}&
\end{flalign}

\begin{remark}
The structure of the function $\mathbb{U}$ given by  (\ref{eq2}) is illustrated in Fig.3.
\end{remark}

\begin{minipage}{11pc}
\setlength{\unitlength}{0.75pt}\begin{picture}(600,220)
\put(30,36){\makebox(0,0)[l]{\footnotesize$0$}}
\put(116,29){\makebox(0,0)[l]{\footnotesize$\sigma$}}
\put(191,29){\makebox(0,0)[l]{\footnotesize$\mathfrak{e}$}}
\put(266,29){\makebox(0,0)[l]{\footnotesize$1$}}
\put(300,29){\makebox(0,0)[l]{\footnotesize$I_{\sigma}^{\mathfrak{e}}$}}
\put(375,29){\makebox(0,0)[l]{\footnotesize$I_{\mathfrak{e}}^{\sigma}$}}
\put(455,29){\makebox(0,0)[l]{\footnotesize$I_{\mathfrak{e},\sigma}$}}

\put(30,69){\makebox(0,0)[l]{\footnotesize$\sigma$}}
\put(30,99){\makebox(0,0)[l]{\footnotesize$\mathfrak{e}$}}
\put(30,129){\makebox(0,0)[l]{\footnotesize$1$}}
\put(25,149){\makebox(0,0)[l]{\footnotesize$I_{\sigma}^{\mathfrak{e}}$}}
\put(25,179){\makebox(0,0)[l]{\footnotesize$I_{\mathfrak{e}}^{\sigma}$}}
\put(20,209){\makebox(0,0)[l]{\footnotesize$I_{\mathfrak{e},\sigma}$}}

\put(76,53){\makebox(0,0)[l]{\footnotesize$0$}}
\put(150,53){\makebox(0,0)[l]{\footnotesize$0$}}
\put(76,85){\makebox(0,0)[l]{\footnotesize$0$}}
\put(131,85){\makebox(0,0)[l]{\footnotesize$\mathbb{U}^{*}(\iota , \kappa )$}}

\put(76,115){\makebox(0,0)[l]{\footnotesize$\iota $}}
\put(76,145){\makebox(0,0)[l]{\footnotesize$0$}}
\put(76,175){\makebox(0,0)[l]{\footnotesize$0$}}
\put(76,205){\makebox(0,0)[l]{\footnotesize$0$}}

\put(131,115){\makebox(0,0)[l]{\footnotesize$\mathbb{U}^{*}(\iota , \kappa )$}}
\put(131,145){\makebox(0,0)[l]{\footnotesize$\mathbb{U}^{*}(\iota , \kappa )$}}
\put(150,175){\makebox(0,0)[l]{\footnotesize$0$}}
\put(150,205){\makebox(0,0)[l]{\footnotesize$0$}}

\put(225,55){\makebox(0,0)[l]{\footnotesize$\kappa $}}
\put(206,85){\makebox(0,0)[l]{\footnotesize$\mathbb{U}^{*}(\iota , \kappa )$}}
\put(206,115){\makebox(0,0)[l]{\footnotesize$\mathbb{U}^{*}(\iota , \kappa )$}}
\put(206,145){\makebox(0,0)[l]{\footnotesize$\mathbb{U}^{*}(\iota , \kappa )$}}
\put(225,175){\makebox(0,0)[l]{\footnotesize$\kappa $}}
\put(225,205){\makebox(0,0)[l]{\footnotesize$\kappa $}}

\put(305,55){\makebox(0,0)[l]{\footnotesize$0$}}
\put(281,85){\makebox(0,0)[l]{\footnotesize$\mathbb{U}^{*}(\iota ,\kappa )$}}
\put(281,115){\makebox(0,0)[l]{\footnotesize$\mathbb{U}^{*}(\iota ,\kappa )$}}
\put(281,145){\makebox(0,0)[l]{\footnotesize$\mathbb{U}^{*}(\iota ,\kappa )$}}
\put(305,175){\makebox(0,0)[l]{\footnotesize$0$}}
\put(305,205){\makebox(0,0)[l]{\footnotesize$0$}}

\put(375,55){\makebox(0,0)[l]{\footnotesize$0$}}
\put(375,85){\makebox(0,0)[l]{\footnotesize$0$}}
\put(375,115){\makebox(0,0)[l]{\footnotesize$\iota $}}
\put(375,145){\makebox(0,0)[l]{\footnotesize$0$}}
\put(375,175){\makebox(0,0)[l]{\footnotesize$0$}}
\put(375,205){\makebox(0,0)[l]{\footnotesize$0$}}

\put(450,55){\makebox(0,0)[l]{\footnotesize$0$}}
\put(450,85){\makebox(0,0)[l]{\footnotesize$0$}}
\put(450,115){\makebox(0,0)[l]{\footnotesize$\iota $}}
\put(450,145){\makebox(0,0)[l]{\footnotesize$0$}}
\put(450,175){\makebox(0,0)[l]{\footnotesize$0$}}
\put(431,205){\makebox(0,0)[l]{\footnotesize$\iota \wedge \kappa \wedge q$}}

\put(45,39){\line(0,1){180}}
\put(120,39){\line(0,1){180}}
\put(195,39){\line(0,1){180}}
\put(270,39){\line(0,1){180}}
\put(345,39){\line(0,1){180}}
\put(420,39){\line(0,1){180}}
\put(495,39){\line(0,1){180}}

\put(270,39){\line(1,0){225}}
\put(270,39){\line(-1,0){225}}

\put(270,71){\line(1,0){225}}
\put(270,71){\line(-1,0){225}}

\put(270,101){\line(1,0){225}}
\put(270,101){\line(-1,0){225}}

\put(270,131){\line(1,0){225}}
\put(270,131){\line(-1,0){225}}

\put(270,161){\line(1,0){225}}
\put(270,161){\line(-1,0){225}}

\put(270,191){\line(1,0){225}}
\put(270,191){\line(-1,0){225}}

\put(270,219){\line(1,0){225}}
\put(270,219){\line(-1,0){225}}

\put(150,0){\emph{Fig.3. The function $\mathbb{U}$ given by (\ref{eq2}).}}
\end{picture}
\end{minipage}\\

In the following, we discuss the function $\mathbb{U}$ given by (\ref{eq2}) with $p\in (\mathfrak{e},1)\cup I_{\mathfrak{e}}^{\sigma}$ and $p\in I_{\sigma}^{\mathfrak{e}}$, respectively.

First, we illustrate that the function $\mathbb{U}$ given by (\ref{eq2}) with $p\in(\mathfrak{e},1)\cup I_{\mathfrak{e}}^{\sigma}$ and $\mathbb{U}^{*}$   can be a uninorm on bounded lattices under some conditions. That is, the dual result of Theorem \ref{th31} is given.

%First, we illustrate that the function $U$ given by (\ref{eq2}) with $p\in(e,1)\cup I_{e}^{b}$   can be a uninorm on bounded lattices under some conditions. That is, the dual result of Theorem \ref{th31} is given.

\begin{theorem}\label{th34}
Let $\sigma \in \mathbb{L_{B}} \setminus\{0,1\}$, $p\in(\mathfrak{e},1)\cup I_{\mathfrak{e}}^{\sigma}$,
$\mathbb{U}^{*}$ be a uninorm on $[\sigma,1]$ with a neutral element $\mathfrak{e}$ and $U_{[\sigma,1]}^{3}$ be a function given by $(\ref{eq2})$. Suppose that  $\iota \wedge \kappa =0$ for all $\iota , \kappa \in I_{e,\sigma}$ with $\iota \neq \kappa $, and $\iota \wedge p =0$ for all $\iota \in I_{e,\sigma}\cap I_{p}$.

$(1)$ Let us assume that $\mathbb{U}^{*}\in \mathcal{U}_{t}$. Then $\mathbb{U}_{[\sigma,1]}^{3}$ is a uninorm on $\mathbb{L_{B}}$ with the neutral element $\mathfrak{e} \in \mathbb{L_{B}} $ if and only if $\iota \parallel \kappa $ for all $\iota \in I_{\mathfrak{e} ,\sigma}\cap I^{p}$ and $\kappa \in I_{\mathfrak{e}}^{\sigma}$.

$(2)$ Let us assume that $I_{\sigma}^{\mathfrak{e}}\cup I_{\mathfrak{e},\sigma}\cup (0,\sigma)\neq\emptyset$. Then $\mathbb{U}_{[\sigma,1]}^{3} $ is a uninorm on $\mathbb{L_{B}}$ with the neutral element $\mathfrak{e} \in \mathbb{L_{B}} $ if and only if $\mathbb{U}^{*}\in \mathcal{U}_{t}$ and $\iota \parallel \kappa $ for all $\iota \in I_{\mathfrak{e},\sigma}\cap I^{p}$ and $\kappa \in I_{\mathfrak{e}}^{\sigma}$.

\end{theorem}

\begin{proof}
 It can be proved with   the proof  of  Theorem \ref{th31} in a similar way.
 \end{proof}\

If we take $e=1$ in Theorem \ref{th34}, then we can obtain the existing result.

\begin{remark}\label{re34}
In Theorem \ref{th34}, if taking $\mathfrak{e}=1$, then $\mathbb{U}^{*}$ be $t$-norm, $\mathbb{U}_{[\sigma,1]}^{3}$ also be $t$-norm, and $I_{\mathfrak{e}}=I_{1}=\varnothing$. In this case, the condition in Theorem \ref{th34} naturally holds.

By the above fact, if taking $\mathfrak{e} =1$ in Theorem \ref{th34}, then we obtain the $t$-conorm $\mathbb{U}_{[\sigma,1]}^{3}:\mathbb{L_{B}}^{2}\rightarrow \mathbb{L_{B}}$ as follows:

$\mathbb{U}_{[\sigma,1]}^{3}(\iota , \kappa )=\begin{cases}
 \mathbb{U}^{*}(\iota ,  \kappa ) &\mbox{if } (\iota , \kappa )\in [\sigma,1]^{2},\\
 \iota \wedge \kappa &\mbox{if } 1\in \{\iota , \kappa \},\\
 0 &\mbox{} otherwise.
\end{cases}$\

Obviously, $\mathbb{U}_{[\sigma,1]}^{3}$ is just the $t$-norm $\mathbb{T}$ in Theorem \ref{th021}.

\end{remark}

\begin{remark}\
In Theorem \ref{th34},  the conditions that $\iota \wedge \kappa =0$ for all $\iota , \kappa \in I_{\mathfrak{e},\sigma}$ with $\iota \neq \kappa $, and $\iota \wedge p =0$ for all $\iota \in I_{\mathfrak{e},\sigma}\cap I_{p}$ can not be omitted, in general.
\end{remark}

\begin{remark}\label{re35}
Let $\mathbb{U}_{[\sigma,1]}^{3}$ be a uninorm in  Theorem \ref{th34}.\\
$(1)$ If $\sigma=0$, then $\mathbb{U}_{[\sigma,1]}^{3}=\mathbb{U}^{*}$.\\
$(2)$ If $\mathbb{U}^{*}\in \mathcal{U}_{t}$, then $\mathbb{U}_{[\sigma,1]}^{3}\in \mathcal{U}_{t}$.\\
$(3)$ $\mathbb{U}_{[\sigma,1]}^{3}\in \mathcal{U}_{min}$ if and only if $\mathbb{U}^{*}\in \mathcal{U}_{min}$.
\end{remark}

At last, we illustrate that the function $U$ given by (\ref{eq2}) with $p\in I_{\sigma}^{\mathfrak{e}}$ and $\mathbb{U}^{*}$  can  be a uninorm on bounded lattices under some conditions. That is, the dual result of Theorem \ref{th33} is given.

\begin{theorem}\label{th36}
Let $\sigma \in \mathbb{L_{B}} \setminus\{0,1\}$, $p\in I_{\sigma}^{\mathfrak{e}}$,
$\mathbb{U}^{*}$ be a uninorm on $[\sigma,1]$ with a neutral element $\mathfrak{e}$ and $\mathbb{U}_{[\sigma,1]}^{4}$ be a function given by $(\ref{eq2})$. Suppose that $\iota \wedge p =0$ for all $\iota \in I_{\mathfrak{e},\sigma}\cap I_{p}$.

$(1)$ Let us assume that $\mathbb{U}^{*}\in \mathcal{U}_{t}$. Then $\mathbb{U}_{[\sigma,1]}^{4}$ is a uninorm on $\mathbb{L_{B}}$ with the neutral element $\mathfrak{e} \in \mathbb{L_{B}} $ if and only if $\iota \parallel \kappa $ for all $\iota \in I_{\mathfrak{e},\sigma}\cap I^{p}$ and $\kappa \in I_{\mathfrak{e}}^{\sigma} $ .

$(2)$ Let us assume that $I_{\sigma}^{\mathfrak{e}}\cup I_{\mathfrak{e},\sigma}\cup (0,\sigma)\neq\emptyset$. Then $\mathbb{U}_{[\sigma,1]}^{4}$ is a uninorm on $\mathbb{L_{B}}$ with the neutral element $\mathfrak{e} \in \mathbb{L_{B}} $ if and only if $\mathbb{U}^{*}\in \mathcal{U}_{t}$ and $\iota \parallel \kappa $ for all $\iota \in I_{\mathfrak{e},\sigma}\cap I^{p}$ and $\kappa \in I_{\mathfrak{e}}^{\sigma} $.

\end{theorem}

\begin{proof}
 It can be proved with the proof  of  Theorem \ref{th33} in a similar way.
 \end{proof}\

If we take $e=1$ in Theorem \ref{th36}, then we can obtain the existing result in the literature.

\begin{remark}\label{re36}
In Theorem \ref{th36}, if taking $\mathfrak{e}=1$, then $\mathbb{U}^{*}$ is a $t$-norm  and $I_{\mathfrak{e}}=I_{1}=\varnothing$. In this case, the condition in Theorem \ref{th36} naturally holds.

By the above fact, if taking $\mathfrak{e}=1$ in Theorem \ref{th36}, then we obtain the $t$-conorm $\mathbb{U}_{[\sigma,1]}^{4}:\mathbb{L_{B}}^{2}\rightarrow \mathbb{L_{B}}$ as follows:

$\mathbb{U}_{[\sigma,1]}^{4}(\iota , \kappa )=\begin{cases}
 \mathbb{U}^{*}(\iota , \kappa ) &\mbox{if } (\iota , \kappa )\in [\sigma,1]^{2},\\
 \iota \wedge \kappa &\mbox{if } 1\in \{\iota , \kappa \},\\
 0 &\mbox{} otherwise.
\end{cases}$\

Obviously, $\mathbb{U}_{[\sigma,1]}^{4}$ is the same as the $t$-norm $\mathbb{T}$ in Theorem \ref{th021}.

\end{remark}

\begin{remark}\
In Theorem \ref{th36},  the condition that $\iota \wedge p =0$ for all $\iota \in I_{\mathfrak{e},\sigma}\cap I_{p}$ can not be omitted, in general.
\end{remark}

\begin{remark}\label{re37}
% Let $\mathbb{U}_{[\sigma,1]}^{4}$ be a uninorm in  Theorem \ref{th36}.\\
$(1)$ If $\sigma=0$, then $\mathbb{U}_{[\sigma,1]}^{4}=\mathbb{U}^{*}$.\\
$(2)$ If $\mathbb{U}^{*}\in \mathcal{U}_{t}$, then $\mathbb{U}_{[\sigma,1]}^{4}\in \mathcal{U}_{t}$.\\
$(3)$ $\mathbb{U}_{[\sigma,1]}^{4}\in \mathcal{U}_{min}$ if and only if $\mathbb{U}^{*}\in \mathcal{U}_{min}$.
\end{remark}

\section{Conclusions}

%In this paper, we have provided some new construction methods of uninorms on bounded lattices $L$ with some conditions on $L$ and the given operators. Moreover, we provide some examples to illustrate these methods for uninorms on bounded lattices. About the above results, we give some remarks as follows.

In this paper, we study the construction methods for uninorms on bounded lattices via functions with the given uninorms and $q\in \mathbb{L_{B}}$ (or $p\in \mathbb{L_{B}}$). Specifically,  we investigate the conditions under which these functions can be uninorms on bounded lattices  when $q\in (0,\mathfrak{e})\cup  I_{\mathfrak{e}}^{\varrho}$ and  $q\in  I_{\varrho}^{\mathfrak{e}}$    (or $p\in (\mathfrak{e},1)\cup I_{\mathfrak{e}}^{\sigma}$ and $p\in I_{\sigma}^{\mathfrak{e}}$), respectively.
%Besides,   $q=0$,  $q\in I_{e,a}$, $q\in [e,a]$ and $q\in (a,1]$ (or $p=1$, $p\in I_{e,b}$,  $p\in [b,e]$ and $p\in [0,b)$)  in \cite{ZX24},  we  have discussed  all the  cases of $q$ (or $p$) in $L$.
Besides the above cases,  Xiu and Zheng \cite{ZX24} discussed how the functions $U$ given by (\ref{eq1}) and (\ref{eq2}) can be a uninorm with $q\in \{0\}\cup I_{\mathfrak{e}, \varrho }\cup [\mathfrak{e},\varrho]\cup (\varrho,1]$ and $p\in \{1\}\cup I_{\mathfrak{e},\sigma}\cup [\sigma,\mathfrak{e}]\cup [0,\sigma)$, respectively.
Up to now,  the functions  given by (\ref{eq1}) and (\ref{eq2}) have been discussed with all cases of $q\in \mathbb{L_{B}}$ and $p\in \mathbb{L_{B}}$, respectively.   Moreover,   our methods generalize
some methods in the literature. See Remarks \ref{re311}, \ref{re333}, \ref{re34} and \ref{re36}.

%Considering construction methods for uninorms on bounded lattices, if the parameters $p$ and $q$ are included in functions,  then we need to choose the appropriate cases of $p$ and $q$ to guarantee that  the functions are uninorms. Moveover, it is necessary and
%interesting to investigate that how the  functions can be  uninorms  with all cases of $q\in L$ or $p\in L$.

Considering the  construction methods for uninorms on bounded lattices, if  $p$ and $q$ are included in functions,  then we usually need to choose the appropriate cases of $p$ and $q$ to guarantee that  the functions are uninorms. Moveover, it is necessary to investigate that how the functions can be  uninorms  with all cases of $q\in \mathbb{L_{B}}$ or $p\in \mathbb{L_{B}}$. In this case, the  construction methods can be studied comprehensively and then provided a novel
perspective to study the constructions of uninorms and other operators.

\end{document}